\documentclass[a4paper, 10pt]{article}
\usepackage{times}
\usepackage[latin1]{inputenc}
\usepackage[T1]{fontenc}
\usepackage{amssymb}
\usepackage{amsmath}
\usepackage{amsthm}
\usepackage{latexsym}
\usepackage[dvips]{graphicx}
\usepackage[left=1.4cm,right=1.4cm,top=1.8cm]{geometry}

\def\N{{\mathbb N}}
\def\R{{\mathbb R}}

\def\e{{\epsilon}}

\def\G{{\Gamma}}

\def\o{{\omega}}

\def\k{{\kappa}}
\def\p{{\prime}}
\def\for{{\,\,\,\forall \,\,\, }}
\def\ex{{\,\,\, \exists \,\,}}

\def\H{ {\mathcal{H} } }
\def\S{ {\mathcal{S} } }

\def\D{{\mathcal{D}}}

\def\div{{\,\mbox{div}\,}}
\def\curl{{\,\mbox{curl}\,}}
\def\pr{\vskip3pt \nl{\bf Proof \,\,}}

\def\id{\,\mbox{id} \,}

\def\+R{+_{_{ \!\! \R}}}
\def\ex{{\mbox{\scriptsize{ex}}}}

\def\tr{\,\mbox{tr}\,}
\def\minus{{\smallsetminus}}

\def\dstar{{\star\star}}

\DeclareMathAlphabet{\mathpzc}{OT1}{pzc}{m}{it}
\def\nl{\vglue0.3truemm\noindent}
\numberwithin{equation}{section}

\begin{document}

\newtheorem{teo}{Theorem}[section]
\newtheorem{pro}[teo]{Proposition}
\newtheorem{lem}[teo]{Lemma}
\newtheorem{defin}[teo]{Definition}
\newtheorem{oss}[teo]{Remark}
\newtheorem{cor}[teo]{Corollary}

\title{On the one fluid limit for vortex sheets}

\author{ 
Fabio Pusateri\\
\vspace{-.2truecm}
{\footnotesize Courant Institute of Mathematical Sciences }
\\{\footnotesize New York University}
\vspace{-.2truecm}
\\{\footnotesize 251 Mercer Street, New York, N.Y. 10012  (USA) }
\vspace{-.2truecm}
\\{\scriptsize pusateri@cims.nyu.edu}
}

\date{}

\maketitle

\begin{abstract}
\nl
We consider the interface problem between two incompressible and inviscid fluids with constant densities 
in the presence of surface tension.
Following the geometric approach of \cite{shatah1,shatah2}
we show that solutions to this problem converge to solutions 
of the free--boundary Euler equations in vacuum as one of the densities goes to zero.
\end{abstract}

\section{Introduction}

\subsection{Description of the problem and main results}
We consider the interface problem between two incompressible and inviscid fluids 
that occupy domains $\Omega^+_t$ and $\Omega^-_t$ in $\R^n$ ($n \geq 2$) at time $t$.
We assume $\Omega_0^+$ is compact and $\R^n = \Omega^+_t \cup \Omega^-_t \cup S_t$ where $S_t := \partial \Omega^\pm_t$.
We let $v_\pm$, $p_\pm$ and $\rho_\pm > 0$ denote respectively the velocity, 
the pressure and the constant density of the fluid occupying the region $\Omega^\pm_t$.
We assume the presence of surface tension on the interface which is argued on physical basis 
to be proportional to the mean curvature $\k_+$ of the hypersurface $S_t$.
The equations of motion are given by\footnote{Here we are introducing the notation
$f = f_+ \chi_{\Omega_t^+} + f_- \chi_{\Omega_t^-}$ for any $f_\pm$ defined on $\Omega_t^\pm$.}
\begin{equation*}
\tag{E}
\label{E}
\left\{
\begin{array}{ll}
\rho (v_t + v \cdot \nabla v) = - \nabla p   &   x \in \R^n \smallsetminus S_t
\\
\\
\nabla \cdot v = 0    &   x \in \R^n \smallsetminus S_t
\\
\\
v (0,x) = v^0 (x) 
											& x \in \R^n \smallsetminus S_0
\end{array}
\right.
\end{equation*}
with corresponding boundary conditions for the interface evolution and pressure's jump given by
\begin{equation}
\tag{BC}
\label{BC}
\left\{
\begin{array}{l}
\partial_t + v_\pm \cdot \nabla  \,\, \mbox{is tangent to} \,\, \bigcup_t S_t \subset \R^{n+1}
\\
\\
p_+ (t,x) - p_- (t,x) = \e^2 \k_+ (t,x) \,\, , \,\,\, x \in S_t  \, .
\end{array}
\right.
\end{equation}
We are interested in analyzing the asymptotic behavior of solutions of the above equations when 
$\rho_- \rightarrow 0$.
Our result is convergence to the solution $(v_+, S_t^\infty)$ of the free--boundary problem for Euler equations 
%
%
\begin{equation*}
\label{E_0}
\tag{$\mbox{E}_0$}
\left\{
\begin{array}{ll}
\rho_+ ( \partial_t v_+ + v_+ \cdot \nabla v_+) = - \nabla p_+   &   x \in \Omega^\infty_t
\\
\\
\nabla \cdot v_+ = 0    &   x \in \Omega^\infty_t
\\
\\
v_+ (0,x) = v^0_+ (x) 
											& x \in \Omega^+_0
\end{array}
\right.
\end{equation*}
with corresponding boundary conditions
\begin{equation*}
\label{BC_0}
\tag{$\mbox{BC}_0$}
\left\{
\begin{array}{l}
\partial_t + v_+ \cdot \nabla  \,\, \mbox{is tangent to} \,\, \bigcup_t S_t^\infty \subset \R^{n+1}
\\
\\
p_+ (t,x) = \e^2 \k^\infty (t,x) \,\, , \,\,\, x \in S_t^\infty 
\end{array}
\right.
\end{equation*}
where $\k^\infty$ denotes the mean curvature of $S_t^\infty := \partial \Omega_t^\infty$.
More precisely we will show the following
\begin{teo}
\label{maintheo}
Let an initial hypersurface\footnote{
The regularity of hypersurfaces in $\R^n$ is intended in the sense of local coordinates:
an hypersurface is $H^s$ for $s > \frac{n}{2}$ if it can be locally represented
as the graph of $H^s$--functions.
} 
$S_0 \in H^{l+1}$ and an initial velocity field $v_0 \in H^l (\R^n \minus S_0)$ 
be given for some $l > \frac{n}{2} + 2$.
Consider any sequence of local in time solutions of (\ref{E})--(\ref{BC}) 
\begin{equation*}
S_t^m \in C ( [0,T], H^{l+1} )
\hskip 8pt , \hskip10pt
v^m \in C ( [0,T], H^l (\Omega_t^m) ) 
\end{equation*}
corresponding to values of the density $\rho^m = \rho_+ \chi_{\Omega_t^+}  +  \rho_-^m \chi_{\Omega_t^-}$
with $\rho^m_- \rightarrow 0$ as $m \rightarrow \infty$.
Then ($v_+^m$, $S_t^m$) converge\footnote{
Convergence is achieved by reducing the problem to the fixed initial domain $\Omega_0$
using Lagrangian coordinate maps. See section \ref{secproof} for details.
}
on a small time interval to the solution
\begin{equation*}
S_t = \partial \Omega_t^\infty \in C ( H^{l^\p + 1} )
\hskip 8pt , \hskip10pt
v_+ \in C ( H^{l^\p} (\Omega_t^\infty) ) \hskip8pt \mbox{for any $l^\p < l$}
\end{equation*}
of (\ref{E_0})--(\ref{BC_0}).
Convergence is in the space
$(S_t, v) \in L^\infty ( H^{l-\frac{1}{2}} ) \times L^\infty ( H^{l^\p} ) $ for any $l^\p < l$.
\end{teo}

\nl
Free--boundary problems for Euler equations have been extensively studied in recent years following the breakthrough
of Wu in \cite{Wu1,Wu2} where local well--posedness in Sobolev spaces is proved in $2$ and $3$ dimensions
for the irrotational gravity water wave problem.
Many works have dealt with the water wave problem also in the general non--zero curl case,
see for instance \cite{Lindblad,ChriLind,shatah1,CoutShko2}.

\nl
For the irrotational vortex sheet problem with surface tension Ambrose \cite{Ambrose} and
more recently Ambrose and Masmoudi \cite{AmbroseMasmoudi2} proved well--posedness respectively in $2$ and $3$ dimensions.
Cheng, Coutand and Shkoller \cite{CoutShko1} proved well--posedness in 3--d for the full problem with rotation
and well--posedness is also obtained (in any dimension) by Shatah and Zeng \cite{shatah3}
for \eqref{E}--\eqref{BC} and other realted fluid surface problems \cite[sec. 6]{shatah3}.

\nl
In absence of surface tension the vortex sheet problem for the free--boundary motion of two fluids is ill--posed
due to the Kelvin--Helmotz instability as shown in \cite{Ebin2}.
Beale, Hou and Lowengrub \cite{Beale} showed how the surface tension regularizes the linearized problem.
In the next section we will show how the Kelvin--Helmotz instability
is very apparent from the infinite--dimensional geometric arguments presented by Shatah and Zeng in \cite{shatah2}.

\nl
We recall that also the free--boundary problem for Euler equations in vacuum (\ref{E_0})--(\ref{BC_0}) with $\e=0$
is known to be ill--posed due to Rayleigh--Taylor instability, see \cite{Ebin1},
which occurs if one does not assume the sign condition
\begin{equation}
\label{RT}
\tag{RT}
- \nabla_{N_+} p_+ (x,t) \geq a > 0   \,\, .
\end{equation}
In \cite{shatah1} it is shown how also the Rayleigh--Taylor instability
is a natural consequence of a geometric calculation
and is related to the sign of an operator appearing in the linearization of the Euler flow.
Motivated by this we are going to show
\begin{pro}
\label{theocurvature}
Let $\G$ be the space of all admissible Lagrangian maps for the interface problem 
\eqref{E}--\eqref{BC} defined in \eqref{defG}
and let 
\begin{equation*}
\G^\star := \left\{ \Phi : \Omega_0^+ \rightarrow \R^n \,\,
										\mbox{\textnormal{volume--preserving homeomorphisms}}
										\right\}
\end{equation*}
be the corresponding space for the water wave problem \eqref{E_0}--\eqref{BC_0}.
Consider a point $u \in \G$ and tangent vectors $v_i \in T_{u} \G$ for $i = 1 \dots 4$,
where $T_{u} \G$ is endowed with the $L^2 (\rho^m dx)$ metric.
If we denote\footnote{
Covariant differentiation on $T_u \G$ (and on $T_u \G^\star$) is defined in section \ref{secgeometryG}.
}
$\bar{\mathcal{R}}$ and $\bar{\mathcal{R}}^\star$ the curvature tensors of $\G$ and $\G^\star$ respectively, then
\begin{equation}
\label{curvatureconv}
{\langle  \bar{\mathcal{R}} (u) (v_1 , v_2) v_3 , v_4  \rangle}_{L^2 (\rho^m dy)}
										\stackrel{m \rightarrow \infty}{\longrightarrow}  
										{\langle \bar{\mathcal{R}}^\star (u_+) (v_{1 +} , v_{2 +} ) v_{3 +} , v_{4 +}  \rangle}_{L^2 (\rho_+ dy)}
\end{equation}
\end{pro}

\nl
In view of the geometric frame work described below and the linearized equation \eqref{linearization},
proposition \ref{theocurvature} can be considered as a first step in showing that solutions of \eqref{E}--\eqref{BC}
converge to solutions of \eqref{E_0}--\eqref{BC_0} with $\e = 0$
when $\e, \rho_- \rightarrow 0$ at the same time\footnote{
We believe that some condition of the form $\rho_- = O( \e^\alpha )$ for some $\alpha > 0$ 
should be needed in this case.}.

\vskip5pt
\nl
Our paper is organized as follows. 
The geometry of $\G$ is presented in section \ref{secgeometry} and
an explanation of the geometric intuition behind the Kelvin--Helmotz and Raileigh--Taylor instabilities  
is given in \ref{seclinearization}. 
Of course we refer to \cite{shatah1,shatah2} for full details about this general geometric approach.
In section \ref{secEE} we state theorems on energy estimates which are independent of $\rho_-$.
Proofs are performed in section \ref{proofVS1} and appendix \ref{proofVS2}.
Section \ref{secproof} is devoted to showing strong convergence of solutions as stated in theorem \ref{maintheo}.
The proof of proposition \ref{theocurvature} is then performed in section \ref{seccurvature}.

\vskip5pt
\nl
During the writing of this manuscript it was brought to the attention of the author that 
Cheng, Coutand and Shkoller \cite{coutshko3} had proved an analogous result to the one stated in theorem \ref{maintheo}.

\subsection{The geometric approach to Euler equations }
\label{secgeometry}
It is well--known that the interface problem between two fluids
has a variational formulation on a subspace of the space of volume--preserving homeomorphisms.
For the water wave problem this was observed for the first time by Arnold in his seminal paper \cite{Arnold66},
where he pointed out that Euler equations for the motion 
of an inviscid and incompressible fluid can be viewed as the geodesic flow 
on the infinite--dimensional manifold of volume--preserving diffeomorphisms.
This point of view has been adopted by several authors in works 
such as \cite{Shnirelman,Brenier,EbinMarsden}
and more recently by Shatah and Zeng in \cite{shatah1,shatah2,shatah3}.

\subsubsection{Lagrangian formulation}
The surface tension parameter $\e$ will be henceforth set to be one.
Multiplying (\ref{E}) by $v$, integrating over $\R^n \minus S_t$, using the boundary condition (\ref{BC})
and the variation of surface area formula, we obtain the conserved energy\footnote{
Notice that the conserved energy does not control 
the $L^2$ norm of $v_-$ in the asymptotic regime $\rho_- \rightarrow 0$.
}
\begin{equation}
\label{energy0}
E = E_0 (S_t, v) = \int_{\R^n \minus S_t} \frac{\rho {|v|}^2 }{2} \, dx + \int_{S_t} \, dS =: 
										\int_{\R^n \minus S_t} \frac{\rho {|v|}^2 }{2} \, dx + S (S_t) \, .
\end{equation}
For $y \in \Omega_0^\pm$ we define $u_\pm (t,y)$ to be the Lagrangian coordinate map associated to
the velocity field $v_\pm$, i.e the solution of the ODE
\begin{equation}
\label{lagrangianmap}
\frac{d x}{dt} = v_\pm (t, x) \hskip8pt, \hskip10pt  x(0,y)= y \hskip 10pt \for y \in \Omega_0^\pm \, ;
\end{equation}
for any vector field $w$ on $\R^n \minus S_t$ we define its material derivative by
\begin{equation*}
\mathbf{D}_t w := w_t + v \cdot \nabla w =  {(w \circ u)}_t \circ u^{-1} \, .
\end{equation*}
In \cite[sec. 2]{shatah2} the authors derive from (\ref{E})--(\ref{BC}) 
the following equation for the physical pressure: 
\begin{equation}
\label{physicalp}
\left\{
\begin{array}{lcl} 
- \Delta p & = & \rho \tr (Dv^2)  
\\
\\
\left. p_\pm \right|_{S_t}  &  =  &
						\mathcal{N}^{-1} \left\{ - \frac{1}{\rho_\mp} \mathcal{N}_\mp  \k_\mp
						- 2 \nabla_{v_+^\top - v_-^\top} v_+^\bot
				    - \Pi_+ ( v_+^\top, v_+^\top ) - \Pi_- ( v_-^\top, v_-^\top ) \right.
				    \\
				    \\
				    &  &  \left. 
				    - \nabla_{N_+} \Delta_+^{-1}  \tr (Dv^2)  - \nabla_{N_-} \Delta_-^{-1}  \tr (Dv^2)
				    \right\}
\end{array}
\right.
\end{equation}
where $\Pi_\pm$ denotes the second fundamental form of the hypersurface $S_t$ (with respect to the outward unit normal vector $N_\pm$ relative to the domain $\Omega_t^\pm$) and $\mathcal{N}$ is given by
\begin{equation}
\label{N}
\mathcal{N} = \frac{\mathcal{N}_+}{\rho_+}  +   \frac{\mathcal{N}_-}{\rho_-}
\end{equation}
with $\mathcal{N}_\pm$ denoting the Dirichlet--to--Neumann operator on the domain $\Omega^\pm_t$.

\nl
From (\ref{lagrangianmap}) we see that in Lagrangian coordinates Euler equations  assume the form
\begin{equation}
\label{Elagrangian}
\rho u_{tt} = - \nabla p \circ u \hskip 20pt u(0) = \id_{\Omega_0}
\end{equation}
with $p$ determined by (\ref{physicalp}).

\nl
Since $v$ is divergence free, $u_\pm$ are volume--preserving maps on $\R^n \minus S_0$.
Moreover $u_+(t, S_0) = u_- (t, S_0)$ even if the restriction to $S_0$ of $u_+$ and $u_-$ do not coincide in general.
This leads to the definition of the space $\G$ of admissible Lagrangian maps for the interface problem:
\begin{eqnarray}
\nonumber
\G & = & \left\{ \Phi = \Phi_+ \chi_{\Omega_0^+} + \Phi_- \chi_{\Omega_0^-} \,\, : \,\,
										\Phi_\pm : \Omega_0^\pm \rightarrow \Phi_\pm (\Omega_0^\pm) \,\, 
										\mbox{is volume--preserving homeo. and}
										\,\,
										\partial \Phi_\pm (\Omega_0^\pm) = \Phi_\pm (\partial \Omega_0^\pm)  \right\} \, .
\label{defG}
\end{eqnarray}
Denoting $S (\Phi) = \int_{\Phi(S_0)} dS$
we can rewrite the energy (\ref{energy0}) in Lagrangian coordinates as
\begin{equation*}
E_0 (u, u_t) = \int_{\R^n \minus S_0} \frac{ \bar{\rho} {|u_t|}^2 }{2} \, dy + S(u)
\end{equation*}
where $\bar{\rho} = \rho \circ u$.
The conservation of the above energy suggests that \eqref{E}--\eqref{BC}
has a Lagrangian action
\begin{equation}
\label{action}
I (u) = \int \int_{\R^n \minus S_0} \frac{\bar{\rho} {|u_t|}^2 }{2} \, dy \, dt \, - \int S(u) \, dt \, .
\end{equation}

\subsubsection{The geometry of $\G$}
\label{secgeometryG}
To derive the Euler--Lagrange equations associated to the action $I$ 
one has to look at the geometry of $\G$ considered as a submanifold  of $L^2 ( \bar{\rho} dy )$
and identify its tangent and normal spaces.
It is easy to see that the tangent space of $\G$ at the point $\Phi$ is given by divergence--free vector fields
with matching normal components in Eulerian coordinates\footnote{We follow the 
convention used in \cite{shatah2} where the Lagrangian description of any vector field 
$X : \Phi (\Omega_0) \rightarrow \R^n$ is denoted by $\bar{X} = X \circ \Phi$.}
\begin{equation}
\nonumber
T_{\Phi} \G = \left\{ \bar{w} : \R^n \minus S_0 \rightarrow \R^n \,\, :
										\nabla \cdot w = 0 \,\, 
										\mbox{ and } \left. w_+^\bot + w_-^\bot \right|_{\Phi(S_0)} = 0 \, , \,\,
										\mbox{where } w  = \bar{w} \circ \Phi^{-1}
										\right\} \, .
\label{defTG}
\end{equation}
while the normal space is
\begin{equation}
\label{TPhiort}
{(T_\Phi \G)}^\bot = \left\{ - (\nabla \psi) \circ \Phi \, : \,  
								\rho_+ \psi_+ \left|_{\Phi(S_0)} \, =  \, \rho_- \psi_- \right|_{\Phi (S_0)} =: \psi^S  \right\} \, .
\end{equation}
A critical path $u(t, \cdot)$ of $I$ satisfies
\begin{equation}
\label{critpath}
\bar{\mathpzc{D}}_t u_t + S^\p (u) = 0
\end{equation}
where $S^\p (u)$ denotes the tangential gradient of $S(u)$ and 
$\bar{\mathpzc{D}}_t$ is the covariant derivative on $\G$ (along $u(t)$).
In order to verify that the Lagrangian map associated to a solution of \eqref{E}--\eqref{BC}
is indeed a critical path of (\ref{action}) we need to
compute $S^\p$ and $\bar{\mathpzc{D}}_t$.


\nl
{\it Computing $\bar{\mathpzc{D}}_t$ and ${II}_{u(t)}$}: 
Given a path $u (t, \cdot) \in \G$ denote $\bar{v} = u_t$ and $S_t = u (t, S_0)$.
For any $\bar{w} (t , \cdot)  \in T_u \G$ we must have
\begin{equation}
\bar{w}_t = \bar{\mathpzc{D}}_t \bar{w}  +  II_{u(t)} (\bar{w}, \bar{v})
\end{equation}
where $II_{u(t)} (\bar{w}, \bar{v}) \in {(T_{u(t)} \G)}^\bot$ denotes the second fundamental form on $T_{u(t)} \G$.
From \eqref{TPhiort} there exists a unique scalar function $p_{v,w}$ defined on $\R^n \minus S_t$ such that 
\begin{equation}
\label{IIfundform}
II_{u(t)} (\bar{w}, \bar{v}) = - \nabla p_{v,w} \circ u \in {\left( T_{u(t)} \G \right)}^\bot
\end{equation}
%
%
%
In \cite{shatah2} it is shown that $p_{v,w}$ is given by\footnote{
Let us point out that in the water wave problem with just one fluid in vacuum we have 
$II_{u(t)}^\star (\bar{w}, \bar{v}) = - \nabla p^\star_{v,w} \circ u \in {\left( T_{u(t)} \G^\star \right)}^\bot$ 
with
\begin{equation}
\label{pwaterwave}
\left\{
\begin{array}{lll}
- \Delta p^\star_{v,w}  & = &   \tr (Dv Dw)    
\\
\\
\left. p^\star_{v,w} \right|_{\partial \Omega_t} & = & 0  \, .
\end{array}
\right.
\end{equation}
}
\begin{equation}
\label{p_vw}
\left\{
\begin{array}{lcl} 
- \Delta p_{v,w} & = & \tr (Dv Dw)
\\
\\
\left. p^\pm_{v,w} \right|_{S_t}  &  =  & \frac{1}{\rho_\pm} p^S_{v,w} =
						- \frac{1}{\rho_\pm} \mathcal{N}^{-1} \left\{ \nabla_{v_+^\top - v_-^\top} w_+^\bot
						+ \nabla_{w_+^\top - w_-^\top} v_+^\bot
				    - \Pi_+ ( v_+^\top, w_+^\top ) \right.
				    \\
				    \\
				    &  &  \left. 
				    - \Pi_- ( v_-^\top, w_-^\top ) 
				    - \nabla_{N_+} \Delta_+^{-1}  \tr (Dv Dw)  - \nabla_{N_-} \Delta_-^{-1}  \tr (Dv Dw)
				    \right\} \, .
\end{array}
\right. 
\end{equation}
Then in Eulerian coordinates we can write
\begin{equation}
\label{covariant}
\mathpzc{D}_t w  :=  \left(  \bar{\mathpzc{D}}_t \bar{w}  \right) \circ u^{-1} = \mathbf{D}_t w + \nabla p_{v,w}  \, .
\end{equation}

\nl
{\it Computing $S^\p (u)$}: For any $\bar{w} \in T_u \G$ the formula for the variation of surface area gives
\begin{equation*}
{ \langle S^\p (u) , \bar{w} \rangle }_{L^2 (\R^n \minus S_0, \rho dy)} = \int_{S_t} \k_+ w_+^\bot \, dS
\end{equation*}
and it is not hard to verify that the unique representation in Eulerian coordinates
of $S^\p (u)$ as a functional acting on $T_u \G$  is

\begin{equation}
\label{Sp}
S^\p (u) = \nabla p_k  \hskip 15pt  \mbox{with}   \hskip 15pt  
													p_k^\pm = \frac{1}{\rho_- \rho_+} \H_\pm \mathcal{N}^{-1} \mathcal{N}_\mp \k_\mp \, .
\end{equation}
From \eqref{physicalp}, \eqref{p_vw} and \eqref{Sp}
we obtain the identity $p = \rho ( p_k + p_{v,v} )$.
Therefore, taking $\bar{w} = u_t$, we see from  \eqref{covariant} and \eqref{Sp}
that a solution of \eqref{critpath} equivalently satisfies
\begin{equation}
\label{Ematerial}
\mathbf{D}_t v  + \nabla p_{v,v} + \nabla p_\k  = 0
\end{equation} 
which is exactly \eqref{Elagrangian} in Eulerian coordinates.

\subsubsection{Linearized equation and instability}
\label{seclinearization}

The Lagrangian formulation discussed above provides a convenient setting to study the linearization of the problem.
Considering variations around the solution $u_t$ of \eqref{critpath} and taking 
a covariant derivative with respect to the variation parameter, 
we obtain the following linearization for $\bar{w} (t, \cdot) \in T_{u(t)} \G$:
\begin{equation}
\label{linearization}
\bar{\mathpzc{D}}^2_t \bar{w} + \bar{\mathpzc{R}} (u) (\bar{u_t}, \bar{w}) u_t + \bar{\mathpzc{D}}^2 S (u) \bar{w} = 0
\end{equation}
where $\bar{\mathpzc{R}}$ denotes the curvature tensor of the manifold $\G$
and $\mathpzc{D}^2 S (u)$ is the projection on $T_{u(t)} \G$ of the second variation of the surface area.
Both of these linear operators acting on $T_u \G$ 
play a central role in the understanding of the problem and in the definition of 
high--order energies based on their leading order terms.
In \cite{shatah1} an explicit but rather complicated formula is given for $\mathpzc{D}^2 S (u)$;
in \cite{shatah1,shatah2} its leading order term $\bar{\mathpzc{A}}$ is singled out and turns out to 
be given\footnote{Both in the one fluid case and the interface problem the leading order term of
$\bar{\mathpzc{D}}^2 S (u)$ has the same form but its Hilbert space representation does not coincide 
due to the different orthogonal splitting of $L^2$ in $T_\Phi \G$ and ${(T_\Phi \G)}^\bot$ in the two settings.
We refer to \cite[pp. 857-858]{shatah2} for the details of the derivation of $\mathpzc{A}$.
}
in Eulerian coordinates by
\begin{equation}
\label{A}
\mathpzc{A}(u) (w) = \nabla f_+ \chi_{\Omega^+}   +    \nabla f_- \chi_{\Omega^-}
\hskip 10pt \mbox{with} \hskip 10pt 
f_\pm = \frac{1}{\rho_+ \rho_-} \H_\pm \mathcal{N}^{-1} \mathcal{N}_\mp (- \Delta_{S_t}) w_\pm^\bot  \, ;
\end{equation}
it is easy to see that $\bar{\mathpzc{A}}$ is a third--order\footnote{
Assuming  $S_t$ is smooth enough.
} 
self--adjoint and positive semi--definite operator with
$ \bar{\mathpzc{A}} (u) ( \bar{w}, \bar{w}) = | \nabla w_\pm^\bot |_{L^2 (S_t)}^2 $.
Further computations performed in \cite[pp 859 - 860]{shatah2}, show that the leading--order 
term $\bar{\mathpzc{R}}_0 (u) (\bar{v})$ of the unbounded sectional curvature operator 
$\bar{\mathpzc{R}}(u) (\bar{v}, \cdot) \bar{v}$
is given in Eulerian coordinates by
\begin{equation*}
\mathpzc{R}_0 (u) (w) = \nabla f_+ \chi_{\Omega^+}   +    \nabla f_- \chi_{\Omega^-}
\hskip 10pt \mbox{with} \hskip 10pt 
f_\pm = \frac{1}{\rho_+ \rho_-} \H_\pm \mathcal{N}^{-1} \mathcal{N}_\mp  \nabla_{v_+^\top - v_-^\top}
\mathcal{N}^{-1} \mathcal{D} \cdot \left( w_\pm^\bot (v_+^\top - v_-^\top) \right) \, .
\end{equation*}
%
%
Noticing that $\bar{\mathpzc{R}}_0 (u)$ is a second--order negative semidefinite differential operator
we immediately see that the linearized Euler equations would be ill--posed
if there had been no surface tension generating the operator $\bar{\mathpzc{A}}$.
This is the so--called {\it Kelvin--Helmotz instability} for the two fluids interface problem.

\nl
We conclude this section recalling that the same geometric setting described above applies
to the problem of Euler equations in vacuum.
The same Lagrangian approach is of course available and the linearized equation is still given by
(\ref{linearization}).
Computations performed in \cite[sec 2.2]{shatah1} 
show how the leading order term of the differential operators involved in the linearization
are given by $\bar{\mathpzc{R}}_0^\star (u)$ and  $\bar{\mathpzc{A}}^\star (u)$ satisfying 
\begin{eqnarray*}
\bar{\mathpzc{R}} (\bar{v}, \bar{w}) & = & \bar{\mathpzc{R}}_0^\star (u) + \, \mbox{bounded operators} 
\\
\bar{\mathpzc{D}}^2 S (u) & = & \bar{\mathpzc{A}}^\star (u) + \, \mbox{second--order differential operators} 
\end{eqnarray*}
and
\begin{equation}
\label{intRT}
\bar{\mathpzc{R}}_0^\star (u) = \int_{S_t}  - \nabla_{N} p^\star_{v,v} {\left| \nabla w^\bot \right|}^2 \, \rho_+ dS
\hskip10pt , \hskip20pt
\bar{\mathpzc{A}}^\star (u) = \int_{S_t} {\left| \nabla w^\bot \right|}^2 \, \rho_+ dS \, .
\end{equation}
Since also in this case $\bar{\mathpzc{A}}^\star (u)$ 
is generated by the presence of surface--tension,
we see that \eqref{linearization} is ill-posed in absence of surface tension if
the {\it Raileigh--Taylor sign condition} \eqref{RT} is not assumed.

\section{Theorems on Energy Estimates}
\label{secEE}

\begin{defin}
\label{defLambda}
Let $\Lambda_0 = \Lambda_0 (S_0, l - \frac{1}{2}, \delta, L)$
for some $l > \frac{n}{2} +1$, $L > 0$ and  $0 < \delta \ll 1$
be the collection of all hypersurfaces $\tilde{S}$
such that a diffeomorphism $F: S_0 \rightarrow \tilde{S} \subset \R^n$ exists with
\begin{equation*}
{| F - \id_{S_0} |}_{H^{l-\frac{1}{2}} (S_0)} < \delta
\end{equation*}
and satisfying a uniform bound on the mean curvature ${|\kappa|}_{ H^{l - \frac{5}{2}} (\tilde{S}) } < L$.
\end{defin}

\nl
In \cite{shatah2} the geometric considerations exposed in section \ref{secgeometry} 
led the authors to define the following energy for \eqref{E}--\eqref{BC}
%
%
%
\begin{defin}
Consider domains  $\Omega_t^\pm$ with $\Omega_t^+$ compact 
and interface $S_t \in H^{l+1}$.
Let $v (t, \cdot) \in H^l (\R^n \minus S_t)$ be any divergence--free vector field 
with $v_+^\bot + v_-^\bot = 0$, define the energy 
\begin{equation}
\label{energyVS}
E (S_t, v(t, \cdot) ) = \frac{1}{2}  \int_{\R^n \minus S_t} {| \mathpzc{A}^{ \frac{l}{3} } v |}^2 \rho \, dx
											+ \frac{1}{2}  \int_{\R^n \minus S_t} {| \mathpzc{A}^{ \frac{l}{3} - \frac{1}{2} } 
											\nabla p_\kappa|}^2 \rho \, dx
											+  {|\o|}^2_{ H^{l-1} (\R^n \minus S_t) } 
\end{equation}
where $\omega$ is the curl of $v$, that is $\omega_i^j = \partial_i v^j - \partial_j v^i$.
\end{defin}

\begin{pro}
\label{proVS1}
Let $l > \frac{n}{2} + 1$, then for $S_t \in \Lambda_0$ with $S_t \in H^{l+1}$ we have
\begin{equation*}
{| \kappa |}^2_{ H^{l-1}(S) } \leq C_0 (1 + E)
\hskip10pt , \hskip20pt
{| v |}^2_{ H^l (\R^n \smallsetminus S) }  \leq   C_0 {(1 + E + E_0)}^2
\end{equation*}
where $C_0$ is some positive constant depending only on $\Lambda_0$ and the initial data 
(in particular it is independent of $\rho_-$).
\end{pro}

\nl
The above proposition is the equivalent of \cite[proposition 4.1]{shatah2}.
The proof of bounds which are independent of the density $\rho_-$ just requires some small modification
of the argument given in \cite{shatah2}.
See section \ref{proofVS1}.

\begin{teo}[\bf{Energy Estimates for (\ref{E}) and (\ref{BC}), \cite{shatah2}}]
\label{teoVS2}
Let $l > \frac{n}{2} + 1$ and a solution to (\ref{E})--(\ref{BC}) be given by 
\begin{equation*}
S_t \in H^{l + 1} \hskip8pt \mbox{and} 
									\hskip8pt v \in C^0_t \left( H^l (\R^n \smallsetminus S_t) \right) \, ,
\end{equation*}
then there exists $L > 0$ and a positive time $t^\star$ independent of $\rho_-$ and depending only on 
${|v(0, \cdot)|}_{H^l (\R^n \smallsetminus S_t)}$, $\Lambda_0$ and $L$,
such that $S_t \in \Lambda_0$ and ${| \kappa |}_{H^{l-1} (S_t) } \leq L$ for all $0 \leq t \leq t^\star$.
Moreover the following energy estimate holds for $0 \leq t \leq t^\star$:
\begin{equation}
\label{energyestimate}
E ( S_t, v(t, \cdot) )  \leq  3 E ( S_0, v(0, \cdot) ) + C_1 + \int_0^t P (E_0, E(S_{t^\p}, v(t^\p, \cdot) )) \, dt^\p
\end{equation}
where $P$ is a polynomial with positive coefficients determined only by $\Lambda_0$ and
the constant $C_1$  depends only on ${|v_0|}_{ H^{l - \frac{3}{2}} (\R^n \smallsetminus S_0) }$ and $\Lambda_0$.
\end{teo}

\nl
The proof of theorem \ref{teoVS2} is essentially the same as in \cite{shatah2} and is postponed to the appendix.
%
%
%
%
\begin{cor} 
\label{unifenergybounds}
Consider a sequence of solutions
\begin{equation*}
S_t^m \in C^0 \left( H^{l+1} \right)  \hskip6pt ,
				\hskip8pt v^m \in C^0 \left( H^l (\R^n \minus S_t^m) \right)
\end{equation*}
solving locally in time the Euler system \eqref{E}--\eqref{BC}
for values of the density $\rho_-^m \rightarrow 0$.
If we denote
\begin{equation*}
E_m (t) := E 	\left( S_t^m, v^m (t, \cdot) \right)
\end{equation*} 
with $E$ given by \eqref{energyVS}, then there exists a positive time $t^\star_0$ 
and a constant $C$ depending only on the set $\Lambda_0$, 
${|v_0|}_{H^l(\R^n \smallsetminus S_0)}$ and ${|v_0|}_{H^{ l-\frac{3}{2} } (\R^n \smallsetminus S_0) }$ such that
\begin{equation}
\label{energybound}
\sup_{ t \in [0, t^\star_0 ] } E_m (t) \leq 2 E(0) + 2 C_1 \hskip 10 pt ,  \hskip 10 pt \for m \in \N \, .
\end{equation}
\end{cor}

\nl
The above corollary gives as a consequence weak convergence of solutions 
of the vortex sheet problem to solutions of the one fluid problem in vacuum.
Standard compactness arguments are going to give the strong convergence stated in theorem \ref{maintheo}. 
See section \ref{secproof} for details.

\nl
For completeness we state here a theorem, proved in \cite{shatah3}, 
based on the above energy estimates and concerning existence of solutions:

\begin{teo}[{\bf Well--posedness for (\ref{E})--(\ref{BC}), \cite{shatah3}}]
\label{teoVS3}
Given an initial surface $S_0 \in H^{l+1}$ and initial velocity $v_0 \in H^l (\R^n \smallsetminus S_0)$
with $l > \frac{n}{2} + 1$, the free interface problem (\ref{E})--(\ref{BC}) has a solution
in the space
\begin{equation*}
S_t \in C^0 \left( H^{l+1} \right)  \hskip8pt , \hskip20pt v \in C^0 \left( H^l (\R^n \smallsetminus S_t) \right)
\end{equation*}
for $t$ in some small interval $[0,T]$ independent of the density $\rho_-$.
Moreover, if $l > 3$ the problem is locally well--posed,
i.e. the solution is unique and depends continuously on the initial data.
\end{teo}

\section{Proof of Proposition \ref{proVS1}}
\label{proofVS1}
Using the definition of $\mathpzc{A}$ in \eqref{A}
we can explicitly write the terms appearing in the energy \eqref{energyVS}
as in \eqref{E_1}, \eqref{E_2} and \eqref{barN} with $k/2$ replaced by $l/3$.
From the properties of $\mathcal{N}_\pm$ and $\mathcal{N}^{-1}$ 
in lemma \ref{lemmaoperator} it follows that there exists a constant $C$ independent of $\rho_-$ such that
\begin{equation*}
{|\k|}^2_{H^{l-1}(S_t)} \leq C {(1+E)}  \hskip7pt,  
\hskip10pt 
{|v_\pm^\bot|}^2_{H^{l-\frac{1}{2}}(S_t)} \leq C {(1 + E )} \, .
\end{equation*}
To estimate $v$ we proceed in three simple steps:
\nl
{\it 1) Estimates on the Lagrangian coordinate map}:
Consider the Lagrangian map $u_-$ associated to $v_-$. 
From lemma \ref{circu} we get
\begin{eqnarray*}
{| u_- (t,\cdot)  -  \id |}_{H^s (\Omega_0^-) } \leq 
								C_1 \int_0^t {| v_- (s, \cdot) |}_{H^s (\Omega_t^-) }
														 {| u_- (s, \cdot) |}_{H^s (\Omega_0^-) }^s \, ds 
\end{eqnarray*}
for any $ 0 \leq s \leq l$ where $C_1 > 0$ only depends on $n$ and $l$.
Now, let $\mu$ be a sufficiently large number to be specified later
depending only on the $H^l$--norm of the initial velocity, define
\begin{equation}
\label{t_0}
t_0 := \sup  \left\{ t \, : \, {|v(t^\p, \cdot)|}_{H^{l} (\R^n \minus S_{t^\p}) } \leq \mu \, \for t^\p \in [0,t]  \right\} \, .
\end{equation}
Since $v$ is assumed to be continuous in time with values in $H^l$, $t_0 > 0$.
The previous inequality and an easy bootstrap argument (or Gronwall's inequality)
show that there exists a positive time $t_1^-$ and a constant $C_2$
depending only on $l,n$,$\mu$ and $\Lambda_0$ such that
\begin{equation}
\label{t_1^-}
{| u_- (t, \cdot) - \id_{\Omega_0^-} | }_{H^l (\Omega_0^-) } \leq C_2 t  \leq \frac{1}{2}
																			\hskip8pt, \hskip10pt \for \, t \in [0, t^\star]
\end{equation}
for $t^\star := \min\{t_0,t_1^-, 1/(2C_2) \}$.
This shows that $u_-$ is an $H^l$--diffeomorphism so that $u_-^{-1} (t, \cdot)$ is 
a well--defined volume preserving map for $x \in \Omega_t^-$ and
for the same range of times we have
\begin{equation*}
{| {(D u_- )}^{-1} | }_{H^s (\Omega_0^-) } \leq 	2   \hskip8pt, \hskip10pt \for \, 0 \leq s \leq l - 1.
\end{equation*}

\vskip5pt
\nl
{\it 2) Decomposition of vector fields and control of the lower norm}:
As it is well-known (and explained in detail in \cite[Appendix B]{shatah1}) any divergence--free
vector field $v : \R^n \minus S_t \rightarrow \R^n$ obeying the condition $v_+^\bot + v_-^\bot = 0$
can be decomposed into two divergence--free components, 
the rotational part $v_r$ responsible for the interior motion
and an irrotational or gradient component $v_{ir} = \nabla g$ responsible for the motion of the boundary $S_t$.
More precisely $g$ is the solution of the elliptic Neumann problem
\begin{equation*}
\left\{
\begin{array}{ll}
\Delta g = 0   &  x \in \R^n \minus S_t
\\
\nabla_{N_\pm} g_\pm = v_\pm^\bot & x \in S_t
\end{array}
\right.
\end{equation*}
and $v_{r} := v - v_{ir}$.
It is observed in \cite{shatah3} that the invariance of Euler equations under the action of 
the group of volume preserving diffeomorphisms leads, via Noether's theorem, to a family of conserved quantities 
which determine the rotational part of the velocity
\begin{equation}
v_r (t,\cdot) = P_r \left( S_t, {(Du^{-1})}^\star v(0, u^{-1}( t, \cdot) )  \right)
\end{equation}
where $P_r (S_t, X)$ denotes the projection of $X : \R^n \minus S_t \rightarrow \R^n$ onto its
rotational (gradient--free) part.
%
%
Therefore we can estimate
\begin{eqnarray*}
{ | v_- (t ,\cdot) | }_{L^2 (\Omega_t^-)} & \leq & { | v_r | }_{L^2 (\Omega_t^-)} 
																						 + { | v_{ir} (t ,\cdot) | }_{L^2 (\Omega_t^-)} 
																						 \leq {| {(Du^{-1})}^\star v(0, u^{-1}( t, \cdot)  |}_{L^2 (\Omega_t^-)}
																						 + { | v_-^\bot (t ,\cdot) | }_{H^\frac{1}{2} (S_t)} 
																						 \\
																						 \\
																						 & \leq & 
																						 {| {(Du)^{-1}} |}_{L^\infty (\Omega_0^-)} 
																						 {|  v(0, \cdot)  |}_{L^2 (\Omega_0^-)} + C E^\frac{1}{2} 
																						 \leq C ( 1 + E^\frac{1}{2})
\end{eqnarray*}
with $C$ depending only on the initial data and $\Lambda_0$.

\nl
{\it 3) Control of ${|v|}_{H^{l}}$}:
To conclude we use the fact\footnote{A more general statement is 
\begin{eqnarray*}
{ |v_\pm| }_{H^l (\Omega_\pm) } & \leq & C ( 1 + {|\k_+|}_{ H^{l - \frac{3}{2}} (S)} ) 
																\left( { |\div v_\pm| }_{H^{l-1} (\Omega_\pm) } 
																+ { |\curl v_\pm| }_{H^{l-1} (\Omega_\pm) } \right.
																+ \left.   {|v_\pm^\bot - \nabla \Delta^{-1} \div v_\pm| }_{
																H^{l-\frac{1}{2}} (S) } + { |v_\pm| }_{L^2 (\Omega_\pm) }
																\right)
\end{eqnarray*}
where the constant $C$ depends only on $\Lambda_0$. 
An essential proof of this can be found in \cite[proposition 4.3]{shatah1}.
} 
that any divergence--free vector field can be controlled by its $\curl$ and normal component:
\begin{eqnarray*}
{ |v| }^2_{H^{l} (\Omega_t^\pm) } & \leq & 
											C  {(1 + {|\k_+|}_{ H^{l- \frac{3}{2}} (S_t) })}^2  \left(  { |\curl v| }^2_{H^{l-1} (\Omega_t^\pm) }
											+  { | v_+ ^\bot | }^2_{ H^{l-\frac{1}{2}} (S_t) }  
											+  { |v| }_{L^2 (\Omega_t^\pm) }^2
											\right)
											\leq  C  {(1 + E + E_0)}^2
\end{eqnarray*}
where the constant $C$ depends only on the initial data and the set $\Lambda_0$ $_\Box$

\section{Proof Theorem \ref{maintheo}}
\label{secproof}
In this section we are going to use the uniform bounds provided by corollary \ref{unifenergybounds} 
combined with the non--linear Eulerian frame work introduced in \cite{shatah1}
to obtain the strong convergence of solutions stated in theorem \ref{maintheo}.

\subsection{Convergence of Lagrangian maps and velocity fields}
\label{secunifestimate}

As a first step we need to estimate the physical pressure.
\begin{lem}
\label{lemp}
Let $v \in H^l$ and $S_t = \partial \Omega_t \in H^{l+1}$ with $l > \frac{n}{2} + 2$
be a given solution of \eqref{E}--\eqref{BC}. 
Then the pressure $p$, determined by \eqref{physicalp}, satisfies
\begin{equation}
\label{estimatep}
{| p_+ |}_{H^{l-\frac{1}{2}} (\Omega_t^+)}  \leq  
												C \left(   {| v_+ |}^2_{H^{l-1} (\Omega_t^+)}  +  {| \kappa_+ |}_{ H^{l - 1} (S_t)}
																	+  \rho_- {| v |}^2_{ H^{l-\frac{1}{2}} (\R^n \minus S_t)} 
																	{| N |}_{ H^{l - 1} (S_t)}  \right)
\end{equation}
and for $\rho_- \ll 1$
\begin{equation}
\label{estimatep-}
{| p_- |}_{H^{l-\frac{1}{2}} (\Omega_t^-)}  \leq  C  \rho_- \left( {| v_- |}^2_{H^{l-1} (\Omega_t^-)} 
																			  + {| \kappa_- |}_{ H^{l - 1} (S_t)}
																			  + {| v |}^2_{ H^{l-\frac{1}{2}} (\R^n \minus S_t)} 
																			  {| N|}_{ H^{l - 1} (S_t)} \right)
\end{equation}
for some constant $C$ depending only on the set of hypersurfaces $\Lambda_0$.
\end{lem}

\pr 
Write $p_\pm = \Delta_\pm^{-1} \Delta p_\pm +  \H_\pm p_\pm^S$ and use lemma \ref{lemmaoperator} 
to get
\begin{eqnarray*}
{| p_\pm |}_{H^{l-\frac{1}{2}} (\Omega_t)} & \leq  & 
 																				C  \left( \rho_\pm {| \tr {(Dv_\pm)}^2 |}_{H^{l-\frac{5}{2}} (\Omega_t)}
																				+ { | p^S_\pm | }_{H^{l-1} (\S_t) }
																				\right)
																				\\
																				\\
																				& \leq & 
																				C \rho_\pm {| v_\pm |}^2_{H^{l-1} (\Omega_t)} 
																				+ C \frac{\rho_-}{\rho_\mp}
																			  {| \kappa_\pm |}_{ H^{l - 1} (S_t)}
																			  + C \rho_- \left( {| N_\pm |}_{H^{l - 2} (S_t)}
																			  	{| v |}^2_{ H^{l-\frac{3}{2}} (\Omega_t)}
																			  \right.
																			  \\
																			  \\
																			  & + & \left. {| v |}_{ H^{l-\frac{3}{2}} (\Omega_t)} 
																			  {| v |}_{ H^{l-\frac{1}{2}} (\Omega_t)} 
																			  {| N_\pm |}_{ H^{l - 1} (S_t)}  \right)  \,\,\,\, _\Box
\end{eqnarray*}

\begin{pro}
\label{lemconvergences}
There exists a sequence $\{ m_k \}$, a time $t^\dstar$ depending only on the initial data 
and an $H^l$--diffeomorphism
$u_+ \in C^0_t \left( [0, t^\dstar] ; H^l (\Omega_0^+) \right)$ 
with $ \partial_t u_+  \in C^0_t \left( [0,t^\dstar] ; H^{l^\p} (\Omega_0^+) \right)$
such that
\begin{eqnarray}
\label{convergenceu}
\lim_{k \rightarrow \infty} u_+^{m_k} = u_+  &  &
							\mbox{in} \hskip 8pt C_t^0 ( \left( [0,t^\dstar] ; H^l (\Omega_0^+) \right) )
\\ \nonumber
\\
\label{convergenceut}
\lim_{k \rightarrow \infty} \partial_t u_+^{m_k} = \partial_t u_+  &  &
							\mbox{in} \hskip 8pt C^0_t \left( [0,t^\dstar] ; H^{l^\p} (\Omega_0^+) \right)
\end{eqnarray}
for any $l^\p < l$.
Moreover if we define 
\begin{equation}
\label{Omegatinfty}
\Omega_t^\infty := u_+ (t, \Omega_0)  
\end{equation}
then there exists
$v_+ \in L^\infty \left( H^l (\Omega_t^\infty) \right) \cap  L^\infty \left( H^{l^\p} (\Omega_t^\infty) \right)$
such that
\begin{eqnarray}
\label{convergencevcircu}
\lim_{k \rightarrow \infty} v_+^{m_k} \circ u_+^{m_k}  =  v_+ \circ u_+ 
							\mbox{in} \hspace{8pt} C_t^0 \left( [0,t^\dstar] ; H^{l^\p} (\Omega_0^+) \right)
\end{eqnarray}
for any $l^\p < l$
and $p_+ \in L^\infty \big( H^{l - \frac{1}{2}} (\Omega_t^+) \big) $
such that
\begin{equation}
\label{convergencepcircu}
\lim_{k \rightarrow \infty} p_+^{m_k} \circ u_+^{m_k}  =  p_+ \circ u_+ 
							\hskip8pt \mbox{weak--star in} \hspace{8pt} L^\infty \left( [0,t^\dstar] ; H^{l-\frac{1}{2}} (\Omega_0^+) \right)
\end{equation}
We will still denote these subsequences by the index $m$.
\end{pro}

\pr
Let us denote 
$X (H^s) = X ( [0, t^\dstar] , H^s (\Omega_0^+) )$ for $X = L^\infty$ or $C_t^0$
and $C$ any positive constant depending only the initial data and the set $\Lambda_0$.
Combining proposition \ref{proVS1} and corollary \ref{unifenergybounds} we see that
\begin{equation*}
{|v^m|}_{L^\infty (H^l (\Omega_t^{+,m})) } \leq C_0 {( 1 + E_m )} \leq C
\end{equation*}
for any $t \leq t^\star_0$. Therefore, arguing as in the proof of proposition \ref{proVS1},
we can find a positive time $t^\dstar \leq t^\star_0$ depending only on $\Lambda_0$ and the initial data, 
such that for any $0 \leq t \leq t^\dstar$
\begin{equation}
\label{time1}
{| u^m_+ (t, \cdot) - \id_{\Omega_0^+} | }_{H^l (\Omega_0^+) } \leq C t^\dstar \leq \frac{1}{2} \, .
\end{equation}
This show that each map $u_+^m$ is an $H^l$--diffeomorphism onto its image and
is uniformly bounded in $L^\infty ( H^l )$ 
by a constant depending only on the initial data and the set $\Lambda_0$.
Then, up to extraction of a subsequence, 
there exist $u_+ \in L^\infty (H^l)$  such that $u_+^m \rightarrow u_+$ weak--star in $L^\infty (H^l)$.
Lemma \ref{circu} and \eqref{time1} imply
\begin{equation*}
{|\partial_t u_+^m|}_{H^l (\Omega_0^+) } \leq {|v_+^m|}_{H^l (\Omega_t^{+,m})} {|u_+^m|}^l_{H^l (\Omega_0^+) } \leq C \, .
\end{equation*}
Again by standard compactness we have, up to extraction,
$\partial_t u_+^m = v^m_+ \circ u_+^m \rightarrow \partial_t u_+ =: \bar{v}_+$  
weak--star in $L^\infty (H^l)$.
Since $u_+^m, u_+ \in W^{1,\infty} (H^l)$, we get $u_+ \in C^0_t (H^l)$
and\footnote{
The standard argument is the following.
Consider an arbitrary subsequence of $\{ u_+^m \}$; the boundedness of $\{ \partial_t u_+^m \}$ implies through
the Ascoli--Arzel\'{a} theorem the existence of a sub-subsequence converging in $C_t^0 (H^l)$ to a limit which
must be $u_+$ (the weak $\star$ limit of the original sequence $\{ u_+^m \}$). 
Therefore $u_+$ is the uniform limit of $\{ u_+^m \}$.
} 
$u_+^m \rightarrow u_+$ in $C^0_t (H^l)$.

\nl
Passing to the limit in \eqref{time1} we see that $u_+$ is also an $H^l$--diffeomorphism.
Thus we can define $v_+$ by $v_+ \circ u_+ =: \bar{v}_+ = \partial_t u_+$.
From Euler equations we have $\partial_t (v^m_+ \circ u_+^m) = -\nabla p_+^m \circ u_+^m$ 
so that lemma \ref{lemp}, lemma \ref{circu}
and corollary \ref{proVS1} together with \eqref{stimaNPI} imply
\begin{eqnarray*}
{| \partial_t (v^m_+ \circ u_+^m) |}_{ H^{l-\frac{3}{2}} {(\Omega_0^+)} } \leq 
						C {| p_+^m |}_{H^{l - \frac{1}{2}} (\Omega_t^{+,m}) }
						\leq  C	 \, .
\end{eqnarray*}
In particular this gives continuity of $v_+^m \circ u_+^m = \partial_t u_+^m$ with values in $H^{l-1} (\Omega_0^+)$.
It also implies the existence of a subsequence (still denoted by the index $m$)
such that $\partial_t (v_+^m \circ u_+^m) \rightarrow \bar{V}_+$ 
weak--star in $L^\infty (H^{l-\frac{3}{2}})$. 
Since $v_+^m \circ u_+^m \rightarrow \bar{v}_+$ in the sense of distributions, $\bar{V}_+ = \partial_t \bar{v}_+$.
Therefore\footnote{
We use the fact that $f \in L^2 (H^{s_1})$ and $f_t \in L^2 (H^{s_2})$
imply $f \in C (H^{ (s_1+s_2)/2 } )$.
} 
$\bar{v}_+ = v_+ \circ u_+ \in C^0_t (H^{l-1})$ and
\begin{equation*}
v_+^m \circ u_+^m \rightarrow v_+ \circ u_+ \hskip8pt
							\mbox{in} \hskip 8pt C_t^0 \left( ( [0,t^\dstar] ; H^{l-1} (\Omega_0^+) ) \right) \, .
\end{equation*}
As $v_+^m \circ u_+^m$ is uniformly bounded in $L^\infty (H^l)$,
by interpolating the Sobolev norms we can improve the above convergence
obtaining \eqref{convergencevcircu} and the equivalent \eqref{convergenceut}.

\nl
Finally, since $ p_+^m \circ u_+^m$ is uniformly bounded in $L^\infty ( H^{l - \frac{1}{2}} )$, 
up to extraction, we have $ p_+^m \circ u_+^m \rightarrow \bar{p}_+ $ weak--star in $L^\infty (H^{ l-\frac{1}{2} })$
and \eqref{convergencepcircu} follows just by defining $p_+ =: \bar{p}_+ \circ u_+^{-1}$
$_\Box$

\subsection{Verification of \eqref{BC_0}}
\label{secmoving}
Using convergence of the Lagrangian maps $u^m_+$ associated to $v^m_+$ established in \eqref{convergenceu},
we defined in \eqref{Omegatinfty} the ``limit domain'' $\Omega_t^\infty$ 
where the evolution of the limit solution is going to take place.
From \eqref{convergenceu} 
and trace estimates we obtain
$ \left. u_+^m \right|_{S_0} \longrightarrow  \left.  u_+   \right|_{S_0} $ in $C^0_t ( H^{ l - \frac{1}{2}} (S_0) )$
so that 
\begin{equation*}
u_+(t, S_0) = \partial u_+(t, \Omega_0^+) =: S_t^\infty 
							\in C^0_t ( H^{ l - \frac{1}{2} } ) \hskip10 pt \mbox{for} \hskip 6pt t \in [0, t^{\dstar}] \, .
\end{equation*}

\begin{pro}
\label{promoving}
The moving boundary condition in \eqref{BC_0} holds for the set of hypersurfaces $S_t^\infty$
with $v_+$ defined by \eqref{convergencevcircu}. 
\end{pro}
\pr
From the definition of Lagrangian maps, \eqref{convergenceut} and \eqref{convergencevcircu}
we have
\begin{equation*}
\partial_t u_+ (t, y) = v_+ ( t , u_+(t,y) ) \hskip8pt \for (t,y) \in [0, t^\dstar] \times \Omega_0^+ \, .
\end{equation*}

\nl
As $u_+ (t, S_0) = S_t^\infty$ for any $t \in [0, t^\dstar]$,
we have that $(t, u_+ (t, \cdot) )$ is a curve on the space--time boundary $ \cup_t S_t^\infty$;
therefore
\begin{equation*}
\partial_t + \partial_t u_+ \cdot \nabla  =  \partial_t + v_+ \circ u_+ \cdot \nabla  
					\,\, \mbox{is tangent to} \,\, \bigcup_t S_t^\infty \subset \R^{n+1} \, .
\end{equation*}
The fact that $u_+$ is a diffeomorphism from $S_0$ to $S_t^\infty$ for any $t \in [0, t^\dstar]$
gives the claim $_\Box$

\begin{lem}
Let $N_+^m (t, \cdot)$ be the outward unit normal and $\kappa_+^m (t, \cdot)$ the mean curvature of $S_t^m$.
Denote by $N^\infty (t,x)$ and $\kappa^\infty (t,x)$ respectively the unit normal and the mean curvature of 
$S_t^\infty$ at the point $x$. Then for any $l^\p < l$
\begin{eqnarray}
\label{limitN}
N_+^m \circ u_+^m \rightarrow N^\infty \circ u_+  \hskip8pt \mbox{in} \,\, C_t^0 ( H^{l^\p} (S_0) )
\hskip 10pt \mbox{and} \hskip15pt
\kappa_+^m \circ u_+^m \rightarrow \kappa^\infty \circ u_+  \hskip8pt \mbox{in} 
																										\,\, C_t^0 ( H^{l^\p - 1} (S_0) )  \, .
\end{eqnarray}
In particular ${|\k^\infty|}_{H^{l^\p - 1} (S_t^\infty)}$ is uniformly bounded
which implies\footnote{This can be proved using local coordinates and estimates for quasi--linear
elliptic equations. Another proof can be found in \cite[proposition A.2]{shatah1}.
}
$S_t^\infty \in H^{l^\p+1}$ as stated in theorem \ref{maintheo}.
\end{lem}

\pr
Since $\k^\infty (t,x) (X,Y) = \tr ( Y \cdot \nabla_X {N^\infty(t,x)} )$ for any $X,Y \in T_x S_t^\infty$, 
it is enough to prove the first statement in \eqref{limitN}.
\nl
We use similar arguments to those in the proof of proposition \ref{lemconvergences}.
By lemma \ref{circu}, \eqref{stimaNPI} and \eqref{energybound} we obtain uniform
bounds on $N_+^m \circ u_+^m$ in $L^\infty (H^l)$;
therefore there exists  $\bar{N}_+ \in L^\infty (H^l)$ such that, up to extraction of a subsequence, 
$N_+^m \circ u_+^m \rightarrow \bar{N}_+ =: A^\infty \circ u_+$ weak--star in $L^\infty (H^l)$.
Identity \eqref{D_tN} and estimate \eqref{stimaNPI} combined with the uniform energy bounds on $\k_+^m$ show that
\begin{equation}
{ \left| \frac{d}{dt} (N_+^m \circ u_+^m)  \right| }_{H^{ l - \frac{3}{2} } (S_0)} \leq 
																								C {|v^m_+|}_{H^l (\Omega_t^m)} {|N_+^m|}_{ H^{l - \frac{3}{2} } (S_0)}
																								\leq C
\end{equation}
with some $C$ uniform in $\Lambda_0$ and $m$.
This in particular implies that $ N_+^m \circ u_+^m $ belongs to $C (H^{l-1} (S_0))$
and that, up to further extraction,
$ \partial_t ( N_+^m \circ u_+^m  ) \rightarrow \partial_t ( A^\infty \circ u_+ ) $
weak--star in $L^\infty H^{l-1} (S_0)$.
As a consequence, $A^\infty \circ u_+ \in C^0_t (H^{l-1} (S_0))$ and
$
N_+^m \circ u_+^m  \rightarrow A^\infty \circ u_+ 
$
in $C_t^0 H^{l^\p} (S_0)$ for any $l^\p < l$.

\nl
To show that $A^\infty (t, \cdot)$ is the outward unit normal $N_+^\infty (t, \cdot)$
to the hypersurface $S_t^\infty$ let $\tau^m \in T_x S_t^m$ be an arbitrary tangent vector. 
Since $u_+^m$ is a diffeomorphism from $S_0$ to $S_t^m$,
there exists a unique tangent vector $\tau_0 \in T_y S_0$ such that $\tau_m = d u_+^m (t,y) \tau_0$,
where $du_+^m (t,y)$  denotes the differential of $u_+^m$ as a map from $S_0$ to $S_t^m$ 
acting on $T_y S_0$ for $y = {(u_+^m )}^{-1} (t,x)$.
Then for any $t \in [0, t^\dstar]$ 
\begin{equation*}
\langle \, N_+^m \left( t, u_+^m(t,y) \right), \, d u_+^m (t,y) \, \tau_0 \, \rangle = 0
\end{equation*}
Letting $m$ go to infinity using \eqref{convergenceu} we obtain
\begin{equation*}
\langle \, A^\infty \left( t, u_+ (t,y) \right), \, d u_+ (t,y) \, \tau_0 \, \rangle = 0
\end{equation*}
Since $\tau_m$, and consequently $\tau_0$, was arbitrarily chosen this
implies that $A^\infty (t,x) \bot T_x S_t^\infty$ for $x = u_+ (t,y)$; 
by the strong convergence established above $A^\infty$ is unitary 
and therefore coincides with $N^\infty (t,x) \, _\Box$

\begin{pro}
The boundary condition \eqref{BC_0} for the pressure is satisfied by the limit solution.
\end{pro}

\pr
For the sequence of solutions $( v^m , S_t^m )$ condition \eqref{BC} holds for every $m \in \N$.
As \eqref{E} is also satisfied for every $m$, the boundary condition for the physical pressure $p_\pm^m$
is the one given in \eqref{physicalp} (where of course every quantity has to be indexed by $m$).
Therefore $(p_+^m - \kappa_+^m ) \circ u_+^m = p_-^m \circ u_+^m$ on $S_0$ and we can use lemma \ref{circu},
\eqref{estimatep-} and trace--estimates to obtain
\begin{eqnarray*}
& & {| (p_+^m - \kappa_+^m ) \circ u_+^m |}_{H^{l-1} (S_0)} \leq 
														C {|p_-^m|}_{H^{l - \frac{1}{2}}(\Omega_t^{-,m})} 
														{|u_+^m|}^{l-\frac{1}{2}}_{H^{l-\frac{1}{2}} (\Omega_0)}
\\
\\
& \leq & C \rho_-^m \left( {| v^m_\pm |}_{H^{ l-1 } (\Omega_t^m)} 
																			  + {| \kappa_-^m |}_{ H^{l - 1} (S_t)}
																			  + {| v^m |}^2_{ H^{l - \frac{1}{2}} (\Omega_t^m)} 
																			  {| N_\pm^m |}_{ H^{l - 1} (S_t^m)} \right) \, . 
\end{eqnarray*}
Since the expression in parentheses above is uniformly bounded by the energies,
letting $m \rightarrow \infty$ and using \eqref{limitN} we get
\begin{equation}
\label{convergencep2}
p_+^m \circ u_+^m \longrightarrow  \k^\infty \circ u_+  \hskip8pt \mbox{in} \,\, C_t^0 H^{l^\p-1} (S_0)
\end{equation}
for any $l^\p < l$.
Using \eqref{convergencepcircu} we conclude that $p_+ (t,x) = \kappa^\infty (t,x)$ 
for any $t \in [0, t^\dstar]$ and $x \in S_t^\infty$ $_\Box$

\subsection{Verification of \eqref{E_0}}

We first need the following estimate:
\begin{lem}  
\label{prop}
Let $p_+$ be given by \eqref{physicalp} then
%
%
\begin{equation}
\label{boundD_tp}
{ | \mathbf{D}_t p_+^m |}_{L^\infty (H^{l-2} (\Omega_t^{+,m}) )} \leq C 	\, .
\end{equation}
for some $C$ uniform in $m$. 
\end{lem}

\pr
In what follows we suppress the use of the index $m$ and let $a \lesssim b$ denote $a \leq C b$
for some constant $C$ independent of $\rho_-$.
Writing $p_+ = \H_+ p_+  +  \Delta^{-1} \tr {(Dv_+)}^2$ we have
\begin{equation}
\mathbf{D}_t p_+  =  \mathbf{D}_t  \H_+ p_+   +  \mathbf{D}_t  \Delta^{-1} \tr {(Dv_+)}^2 
											= \H_+ \mathbf{D}_t p_+   +   \Delta^{-1}  \mathbf{D}_t  \tr {(Dv_+)}^2
											+ R
											:= \mbox{(I)} + \mbox{(II)} + R
\label{dtpcircu}											
\end{equation}
where the remainder is given by the sum of the two commutators
\begin{equation}
\label{Rm}
R = R^1 + R^2 :=   \left[ \mathbf{D}_t , \H_+  \right] p_+
						+ \left[ \mathbf{D}_t , \Delta^{-1}  \right] \tr {(Dv_+)}^2   \, .
\end{equation}

\nl
We show that every term is bounded in $H^{l-2}$ or better
by the quantities ${|v|}_{H^l}$, ${|p|}_{H^{l-\frac{1}{2}}}$, ${|\k|}_{H^{l-1}}$ and ${|N|}_{H^l}$
which are already known to be bounded uniformly in time by the energies independently of $\rho_-$.

\nl
{\it Estimate of} (I): This is the highest order term in \eqref{dtpcircu}.
Denoting $P := \left. \mathcal{N} p_+ \right|_{S_t}$ 
we have
\begin{eqnarray*}
\mbox{(I)}   =  \H_+ \mathbf{D}_t \mathcal{N}^{-1} P 
													 =  \H_+ \mathcal{N}^{-1} \mathbf{D}_t  P
													 +  \H_+ R^3 P
							\hskip15pt \mbox{with}  \hskip 10pt  R^3 := \left[ \mathcal{N}^{-1}, \mathbf{D}_t \right]  \, .
\end{eqnarray*}
Observe that $R^3 = \mathcal{N}^{-1} \left[ \mathcal{N}, \mathbf{D}_t \right]  \mathcal{N}^{-1}$
so that \eqref{estimateDH}, \eqref{estimateN-1} and \eqref{commest1} give
\begin{eqnarray*}
{ | \H_+ R^3 P  | }_{H^{l-\frac{1}{2}} ( \Omega_t^+ ) } 
								  & \lesssim &  {\left|  \left[ \mathcal{N}^{-1}, \mathbf{D}_t \right] P  \right| }_{H^{l-1} (S_t) }
									\lesssim   \rho_-  {| v |}_{ H^l (\Omega_t)}    {|  P  | }_{H^{l - 2} (S_t) }
									\\
									\\
									& \lesssim &    {| v |}_{ H^l (\Omega_t)}
							\left( {| \kappa |}_{ H^{l - 1} (S_t)}
							+ {| v |}^2_{ H^{l - \frac{1}{2} } (\Omega_t)}  {| N |}_{ H^{ l-1 } (S_t)} \right) \, .
\end{eqnarray*}
Using again \eqref{estimateDH} and \eqref{estimateN-1} we obtain
\begin{equation}
\label{D_tP}
{|  \H_+ \mathcal{N}^{-1} \mathbf{D}_t  P   |}_{H^{l-2} (\Omega_t^+)} \leq 
						C  \rho_- {|  \mathbf{D}_t P   |}_{H^{l - \frac{7}{2}} (S_t)} \, .
\end{equation}
Now $\mathbf{D}_t P$  contains four different terms to be estimated.
The term involving the mean curvature is estimated by \eqref{commest1} and \eqref{D_tkappa}:
\begin{eqnarray*}
& & {| \mathbf{D}_t \frac{1}{\rho_-} \mathcal{N}_- \kappa_+ |}_{ H^{l - \frac{7}{2}} (S_t) }  \lesssim
						\frac{1}{\rho_-} \left( {| [ \mathbf{D}_t ,\mathcal{N}_- ] \kappa_+ |}_{ H^{l - \frac{7}{2}} (S_t) } 
						+ 
						{ | \mathbf{D}_t \kappa_+ |}_{ H^{l - \frac{5}{2}} (S_t) }
						\right)
						\\
						\\
						& \lesssim & \frac{C}{\rho_-} \left(  {| v |}_{ H^l (\Omega_t)} { | \kappa_+ |}_{ H^{l - \frac{5}{2}} (S_t) }
						+ { | v_+ |}_{ H^l (\Omega_t^+) } { | N_+ |}_{ H^{l - \frac{5}{2}} (S_t) }
						+ { | \k_+ |}_{ H^{l - \frac{5}{2}} (S_t) } { | v_+ |}_{ H^{l - 1} (\Omega_t^+) } 
						\right) \, .
\end{eqnarray*}
Notice that the presence of $\rho_-$ in the denominator in this last estimate is compensated by 
the factor $\rho_-$ in \eqref{D_tP} so that the bounds remain uniform.
For the terms involving $\tr {(Dv)}^2$  we use \eqref{D_tN}, \eqref{D_tDelta-1}  and
the identities $ \mathbf{D}_t \nabla f = \nabla \mathbf{D}_t f - {(Dv)}^\star \nabla f$ and\footnote{
This identity follows from $ \mathbf{D}_t D v = D \mathbf{D}_t v - {(D v)}^2$
together with Euler equations $\rho \mathbf{D}_t v = - \nabla p$.
} 
$\mathbf{D}_t \tr {(Dv)}^2 = - 2 \tr [ {(Dv)}^3  - 2 \rho_+ D^2p \cdot D v]$
to estimate
\begin{eqnarray*}
& & {| \mathbf{D}_t \nabla_{N_\pm} \Delta_\pm^{-1} \tr {(Dv_\pm)}^2|}_{ H^{l - 2} (S_t) }  
\lesssim
						{| \mathbf{D}_t  N_\pm |}_{ H^{l - 2} (S_t) } {| v |}^2_{ H^{l - \frac{3}{2} } (\Omega_t) } 
						+
						{| N_\pm |}_{ H^{l - 2} (S_t) }
						{| \mathbf{D}_t  \nabla \Delta_\pm^{-1} \tr {(Dv_\pm)}^2|}_{ H^{l - \frac{3}{2}} (\Omega_t) }
\\
\\
&  \lesssim &  {| N_\pm |}_{ H^{l - 2} (S_t) } {| v |}_{ H^{l - \frac{1}{2} } (\Omega_t) } 
					 {| v |}^2_{ H^{l - \frac{3}{2}} (\Omega_t) }
				+  {| N_\pm |}_{ H^{l - 2} (S_t) } \left( 
			     {| [\mathbf{D}_t, \Delta_\pm^{-1}] \tr {(Dv_\pm)}^2|}_{ H^{l - \frac{1}{2} } (\Omega_t) }
			  +  {| \mathbf{D}_t \tr {(Dv_\pm)}^2|}_{ H^{l - \frac{5}{2}} (\Omega_t) }  
			 \right.
\\
\\
&  +  & \left. {| {(Dv_\pm )}^\star \nabla \Delta_\pm^{-1} \tr {(Dv_\pm)}^2|}_{ H^{l - \frac{3}{2}} (\Omega_t) }  
\right)
\\
\\
& \lesssim &	{| N_\pm |}_{ H^{l - 2} (S_t) } {| v |}^3_{ H^{l-\frac{1}{2}} (\Omega_t) }
       +  {| N_\pm |}_{ H^{l - 2} (S_t) } \left( 
			    {| v_\pm |}_{ H^{l} (\Omega_t) }   {| v_\pm |}^2_{ H^{l - 1} (\Omega_t) }
			 +  {| v_\pm |}^3_{ H^{l - 1} (\Omega_t) }
			 +  {| p_\pm |}_{ H^{l - \frac{1}{2}} (\Omega_t) } {| v |}_{ H^{l - 1} (\Omega_t) }
			 \right.
\\
\\
& + &  \left. 
			{| v_\pm |}_{ H^{l - \frac{1}{2}} (\Omega_t) }   {| v_\pm |}^2_{ H^{l - 1} (\Omega_t) }
			\right)  \, .
\end{eqnarray*}
Analogously, using \eqref{D_tPI} the terms $\mathbf{D}_t \Pi_\pm (v_\pm^\top, v_\pm^\top)$ 
and $\mathbf{D}_t ( v_\pm^\top \nabla v_\pm^\bot )$ can be bounded uniformly in 
$H^{l-\frac{5}{2}} (S_t)$ and $H^{l-3}(S_t)$ respectively.

\nl
{\it Estimate of $(II)$}: 
By the same formula used above to express $\mathbf{D}_t  \tr {(Dv_+)}^2$ we get
\begin{eqnarray*}
{| \Delta^{-1}  \mathbf{D}_t  \tr {(Dv_+)}^2 |}_{H^{l - \frac{1}{2} } ( \Omega_t^+ )} 
					& \lesssim & 
					{| {(Dv_+)}^3 |}_{H^{ l - \frac{5}{2} } ( \Omega_t^+ )} +
					\rho_+ {|  D^2 p_+ \cdot D v_+ |}_{H^{ l-\frac{5}{2} } ( \Omega_t^+ )} 
					\\
					\\
					& \lesssim &  \left(  {| v_+ |}^3_{H^{l-1} ( \Omega_t^+ )}
					+  {| p_+|}_{H^{l-\frac{1}{2}} ( \Omega_t^+ )} 
					{| v_+|}_{H^{l-1} ( \Omega_t^+ )}  \right) \, .
\end{eqnarray*}

\nl
{\it Estimate of }$R$: Commutators $R^1$ and $R^2$ 
are estimated directly by \eqref{D_tH} and \eqref{D_tDelta-1}:
\begin{eqnarray*}
{\left|  [ \mathbf{D}_t , \H_+  ] p_+  \right|}_{H^{l-2} ( \Omega_t^{+} )} 
												\lesssim   { | v_+ | }_{H^l ( \Omega_t^{+} )} 
												{ | p_+ | }_{H^{l - \frac{5}{2} } ( \Omega_t^{+} )}
\hskip10pt, \hskip15pt												
{ \left|  [ \mathbf{D}_t , \Delta^{-1}  ]  \tr {(Dv_+)}^2  \right|}_{H^l ( \Omega_t^{+} )} 
												\lesssim   { | v_+ | }^3_{H^l ( \Omega_t^{+} )}
\end{eqnarray*}
where as usual the constant $C$ is independent of $\rho_-$ $_\Box$

\begin{pro}
\label{proE_0}
Let $v_+$ and $u_+$ be given as in proposition \ref{lemconvergences} then
\begin{equation*}
\frac{d}{dt} (v_+^m \circ u_+^m) \longrightarrow \frac{d}{dt} (v_+ \circ u_+)  
								\hskip8pt \mbox{in} \hskip5pt C_t^0 (H^{l^\p - \frac{3}{2} } (\Omega_0^+))  \, .
\end{equation*}
and $v_+$ satisfies Euler equations  \eqref{E_0}.
\end{pro}

\pr
\eqref{convergencevcircu} and the uniform bounds on $p_+^m$ 
establish the above convergence weak--star in $L^\infty (H^{l- \frac{3}{2}} (\Omega_0^+))$.
Since $\partial_t^2 (v_+^m \circ u_+^m) = \mathbf{D}_t \nabla p_+^m \circ u_+^m
																 =  \nabla \mathbf{D}_t p_+^m \circ u_+^m	
																 -  {(Dv_+^m)}^\star \nabla p_+^m \circ u_+^m $
the bound given in \eqref{boundD_tp} implies
$\partial_t^2 (v_+^m \circ u_+^m) \in L^\infty (H^{l-3} (\Omega_0^+))$
and the desired strong convergence follows through the usual arguments.
\nl
From \eqref{convergencepcircu}  and  \eqref{convergencep2}
we know that $p_+^m \circ u_+^m \rightarrow p_+ \circ u_+$ strongly in $C^0_t H^{l^\p - \frac{1}{2}}$
and therefore
$\nabla p_+^m \circ u_+^m = \nabla( p_+^m \circ u_+^m ) {(\nabla u_+^m)}^{-1} \rightarrow
\nabla p_+ \circ u_+$  in $C^0_t H^{l-2} (\Omega_0^+)$.
Since Euler equations  in Lagrangian coordinates are
$\partial_t ( v_+^m \circ u_+^m )  =  - \nabla p_+^m \circ u_+^m$
we can take the limit in $L^\infty (H^{l - 2} (\Omega_0^+))$
obtaining that $v_+$ satisfies Euler equations  in Lagrangian coordinates too, i.e.
\begin{equation*}
\frac{d}{dt} v_+(t,  u_+(t,y) ) =  - \nabla p_+ (t, u_+(t,y) )  \hskip8pt \for (t,y) \in [0, t^\dstar] \times \Omega_0^+  \, .
\end{equation*}
Finally from \eqref{convergenceu} and \eqref{convergencevcircu} we have 
$\nabla (v_+^m \circ u_+^m) \rightarrow  \nabla (v_+ \circ u_+)$ in $C_t^0 (H^{l-2} (\Omega_0^+))$ so that
\begin{eqnarray*}
0 \equiv \nabla \cdot v_+^m \circ u_+^m
			=   \tr( \nabla v_+^m  \circ u_+^m )
			= \tr \left( \nabla (v_+^m  \circ u_+^m) {(\nabla u_+^m)}^{-1} \right)
			\stackrel{m \rightarrow \infty}{\longrightarrow}
			\tr \left( \nabla (\bar{v}_+ \circ u_+) {(\nabla u_+)}^{-1} \right)
			= \nabla \cdot v_+ \circ u_+
\end{eqnarray*}
which implies $\nabla \cdot v_+ = 0$ pointwise in $\Omega_t^\infty$ for any $t \in [0, t^\dstar]$ $_\Box$

\vskip6pt
\nl
The proof of theorem \ref{maintheo} is completed $_\blacksquare$

\section{Proof of Proposition \ref{theocurvature}}

\label{seccurvature}
Let $\G$ be the infinite--dimensional manifold (\ref{defG}) and $\bar{\mathcal{R}}$ its curvature tensor
induced by the covariant differentiation defined in section \ref{secgeometry}.
Consider a map $u (t): \Omega_0 \rightarrow \Omega_t$ in $\G$.
Let $\bar{\mathcal{R}}^m$ denote the sectional curvature of $\G$ at the point $u$
as an operator acting on $T_{u} \G$ endowed with the $L^2 (\rho^m  dy)$ metric  
and depending on some $\bar{v} \in T_{u} \G$ (and of course on $u$).
We assume $v$ and the hypersurfaces $S_t$ to be sufficiently smooth 
and single out the leading order term of $\bar{\mathcal{R}}^m$ analyzing its behavior 
as $m$ goes to infinity (or equivalently as the density $\rho_-$ vanishes).

\nl
In view of the geometrical frame work discussed in section \ref{secgeometry},
and in particular in \ref{seclinearization},
$\bar{\mathcal{R}}^m$ can be considered as a measurement 
of the instability occurring in the linearized Euler equations 
in case surface tension were not present.

\nl
Let $\bar{w}$ be any vector in $\in T_{u} \G$. We assume that $w$ is uniformly bounded in $H^l (\R^n \minus S_t)$
for some large enough $l$ and compute the sectional curvature in the direction of $\bar{v}$ and $\bar{w}$.
Using a well--known formula from Riemannian geometry 
together with (\ref{IIfundform}) we have
\begin{eqnarray*}
\bar{\mathcal{R}}^m & = & 
  {\langle \bar{\mathcal{R}} (u)( \bar{v}, \bar{w} ) \bar{v} \, , \, \bar{w} \rangle}_{L^2 (\rho^m  dx)}
= {\langle II_{u} ( \bar{v}, \bar{v} ) \, , \, II_{u} ( \bar{w}, \bar{w} )  \rangle}_{L^2 (\rho^m  dx)} -  
{ \left\| II_{u} (\bar{v}, \bar{w}) \right\| }_{L^2 (\rho^m  dx)}^2
\\
\\
& = & \int_{ \R^n \smallsetminus S_t} \nabla p_{v,v} \nabla p_{w,w} \, \rho^m \, dx   
					-   \int_{\R^n \smallsetminus S_t} {|\nabla p_{v,w} |}^2 \, \rho^m \, dx \, .
\end{eqnarray*}

\nl
Again we suppress the use of the index $m$.
Using the divergence theorem the first integral can be written as
\begin{eqnarray*}
& & \int_{ \R^n \smallsetminus S_t} \nabla p_{v,v} \nabla p_{w,w} \rho \, dx  
				  =  \int_{S_t} p^S_{v,v} \left(  \nabla_{N_+}  p^+_{w,w}  + \nabla_{N_-} p^-_{w,w}  \right) 
				  - \int_{\R^n \smallsetminus S_t}  p_{v,v} \Delta p_{w,w} \rho \, dx 
				  \\
				  \\
					& = & \int_{S_t} p^S_{v,v} \left\{  - 2 \nabla_{w_+^\top - w_-^\top} w_+^\bot
				  + \Pi_+ ( w_+^\top, w_+^\top ) + \Pi_- ( w_-^\top, w_-^\top )
				  \right\} \, dS
					+ \int_{\R^n \smallsetminus S_t} p_{v,v} \tr{(Dw)}^2   \, \rho \, dx  
\end{eqnarray*}
having used
$ \nabla_{N_+} p^+_{w,w}  + \nabla_{N_-} p^-_{w,w} = \mathcal{N} p^S_{v,v}  + 
\nabla_{N_+} \Delta_+^{-1} \Delta p^+_{w,w}  + \nabla_{N_-} \Delta_-^{-1} \Delta p^-_{w,w}$
and \eqref{p_vw} with $v=w$.
Since $tr{(Dw)}^2 = \partial_i w^k \partial_k w^i = \partial_i (w^k \partial_k w^i)$ 
we can use twice again the divergence theorem obtaining
%
%
%
%
\begin{eqnarray}
\nonumber
\int_{ \R^n \smallsetminus S_t} \nabla p_{v,v} \nabla p_{w,w} \rho \, dx  
				  & = & \int_{\R^n \smallsetminus S_t} D^2 p_{v,v} (w,w) \, \rho \, dx  
				  + \int_{S_t} p^S_{v,v} \left\{  - 2 \nabla_{w_+^\top - w_-^\top} w_+^\bot 
				  + \Pi_+ ( w_+^\top, w_+^\top ) + \Pi_- ( w_-^\top, w_-^\top )
				  \right.
				  \\
				  \nonumber
				  \\
				  \label{sectional1}
				  & + & \left. 
				  \nabla_{w_+} w_+ \cdot	N_+	 +  \nabla_{w_-} w_- \cdot	N_-		\right\}  \, dS
					- \int_{S_t} \rho_+ w_+^\bot \nabla_{w_+} p^+_{v,v}  +  \rho_- w_-^\bot \nabla_{w_-} p^-_{v,v} \, dS 
					\, .
\end{eqnarray}
To estimate the terms containing $p^S_{v,w}$, which is the inverse image through $\mathcal{N}$
of a mean zero function on $S_t$,  we use lemma \ref{lemmaoperator}. 
For any $f \in L^1 (S_t)$, \eqref{estimateN-1} yields
\begin{eqnarray*}
\nonumber
\left|  \int_{S_t}  p^S_{v,w}  f  \, dS \right|  \leq  C {| p^S_{v,w} |}_{H^s (S_t)} {|f|}_{L_1 (S_t)}
								 \leq   C \rho_- {| \mathcal{N} p^S_{v,w} |}_{H^{s-1} (S^t) }  {|f|}_{L_1 (S^t) }  
\end{eqnarray*}
whenever $s > \frac{n-1}{2}$, with $C$ uniform in $S_t \in \Lambda_0$.
Since ${| \mathcal{N} p^S_{v,w} |}_{H^{s-1} (S_t) }$ 
is uniformly bounded for smooth enough and bounded $v$ and $w$,
we can easily estimate several terms in \eqref{sectional1}: 
\begin{eqnarray*}
\left|  \int_{S_t} p^S_{v,v} \nabla_{w_\pm^\top} w_+^\bot \, dS  \right|
\hskip5pt, \hskip7pt
\left|  \int_{S_t} p^S_{v,v} \Pi_\pm ( w_\pm^\top, w_\pm^\top ) \, dS  \right|  
\hskip5pt, \hskip7pt
\left|  \int_{S_t} p^S_{v,v} \nabla_{w_\pm} w_\pm \cdot	N_\pm   \, dS  \right|
							\leq
							C \rho_- {|w|}_{H^\frac{3}{2} (\R^n \smallsetminus S_t)}^2  
\end{eqnarray*}
where $C$ is some uniform constant depending on $v$ and the mean curvature of $S_t$.
These bounds imply
\begin{eqnarray*}
&  & \lim_{\rho_- \rightarrow 0}
							\int_{ \R^n \smallsetminus S_t} \nabla p_{v,v} \nabla p_{w,w} \rho \, dx  
							=
							\int_{ \Omega_t^{+} } D^2 p_{v_+, v_+} (w_+, w_+) \, \rho_+ \, dx 
							-  
							\lim_{\rho_- \rightarrow 0}
   						\int_{S_t} \rho_+ w_+^{\bot} \nabla_{ w_+ } p^+_{v, v}  +  
   						\rho_- w_-^{\bot} \nabla_{w_-} p^-_{v, v} \, dS
   						\, .
\end{eqnarray*}
Next we look at the contribution of ${\|  {II}_u (\bar{v}, \bar{w})  \|}^2$,
use the decomposition  $f_\pm = \H_\pm f + \Delta^{-1}_\pm \Delta f$ applied to $p_{v,w}$ 
and observe that $\nabla \H_\pm \bot \nabla \Delta_\pm^{-1} \Delta $ to obtain
\begin{eqnarray*}
\int_{\R^n \minus S_t} {| \nabla p_{v,w}|}^2 \rho \, dx 
					=
					\int_{S_t} p^S_{v,w} \mathcal{N} p^S_{v,w}
					+ \int_{\Omega^+_t}  {| \nabla \Delta^{-1}_+ \tr( Dv Dw) |}^2 \, \rho_+ \, dx
					+ \int_{\Omega^-_t}  {| \nabla \Delta^{-1}_- \tr (Dv Dw) |}^2  \, \rho_- \, dx \, .
\end{eqnarray*}
In \cite{shatah2} it is shown how the leading order term of the 
sectional curvature comes from the contribution
of the surface integral in the above expression and is a second order negative semi--definite operator.
But since $\mathcal{N} p^S_{v,w}$ is independent of $\rho_-$,
by the same argument performed above 
the boundary integral vanishes as $\rho_- \rightarrow 0$.
Therefore
\begin{eqnarray}
\nonumber
\lim_{m \rightarrow \infty}  \bar{\mathcal{R}}^m & = &   
					\int_{ \Omega_t^{+} } D^2 p_{v, v} (w, w) \, \rho_+
					-  {|  \nabla \Delta^{-1}_+ \tr( Dv Dw ) |}^2  \, \rho_+ \, dx 
\\ \nonumber
\\		
\label{sectional2}
& - & \lim_{\rho_- \rightarrow 0}  
			\int_{S_t} \rho_+ w_+^{\bot} \nabla_{w_+} p^+_{v, v} + 
								\rho_-  w_-^{\bot} \nabla_{w_-} p^-_{v, v} \, dS  
								\, .
\end{eqnarray}
By splitting $w$ into normal and tangential components on the boundary 
the surface integral in (\ref{sectional2}) is
\begin{eqnarray}
\label{sectional3}
\int_{S_t} \rho_+ w_+^\bot \nabla_{w_+} p^+_{v,v} + \rho_- w_-^\bot \nabla_{w_-} p^-_{v,v} \, dS  & = &
\int_{S_t} \rho_+ {|w_+^\bot|}^2 \nabla_{N_+} p^+_{v,v} + \rho_- {|w_-^\bot|}^2 \nabla_{N_-} p^-_{v,v}
\\ \nonumber
\\ 
\label{sectional4}
& + & \int_{S_t} \rho_+ w_+^\bot \nabla_{w_+}^\top p^+_{v,v} + \rho_- w_-^\bot \nabla_{w_-}^\top p^-_{v,v} \, dS \, .
\end{eqnarray}
Writing $\rho_\pm p^\pm_{v,v} = \H_\pm \left(p^S_{v,v} \right) - \rho_\pm \Delta_\pm^{-1} \tr{(Dv)}^2$, 
by the usual estimate for $p^S_{v,v}$
the right--hand side of (\ref{sectional3}) gives the contribution
\begin{eqnarray}
\nonumber
& &  \int_{S_t} {|w_+^{\bot}|}^2 \mathcal{N}_+ p^S_{v, v}  + {|w_-^{\bot}|}^2  \mathcal{N}_- p^S_{v, v}
			+  \rho_+   {|w_+^{\bot}|}^2  \nabla_{N_+} \Delta^{-1}_+ \tr{(Dv)}^2 + 
			\rho_-  {|w_-^{\bot}|}^2  \nabla_{N_-} \Delta^{-1}_- \tr{(Dv)}^2 \, dS
\\ 
\nonumber
\\ 
\label{sectional5}
		&  &  \stackrel{\rho_- \rightarrow 0}{\longrightarrow}    \int_{S_t}  \rho_+ \nabla_{N_+} \Delta^{-1}_+ 
																						\tr{(Dv)}^2 {|w_+^{\bot}|}^2 \, dS \, .
\end{eqnarray}
Since $\rho_+ p_{v,v}^+ = \rho_- p_{v,v}^- = p^S_{v,v}$ on $S_t$ and we are considering only tangential derivatives,
the contribution of the term in \eqref{sectional4} is
\begin{equation}
\left| \int_{S_t} w_+^\bot \nabla_{w_+^\top + w_-^\top} p^S_{v,v} \, dS  \right| 
\leq
C {| w |}^2_{H^1(\R^n \smallsetminus S_t)} {|p^S_{v,v}|}_{H^{s_1} (S_t) }
\leq  C \rho_- {|w|}^2_{H^1 (\R^n \smallsetminus S_t) } \stackrel{\rho_- \rightarrow 0}{\longrightarrow} 0 \, .
\label{sectional6}
\end{equation}
Gathering \eqref{sectional2}, \eqref{sectional3}, \eqref{sectional5} and \eqref{sectional6} we get
\begin{eqnarray}
\label{sectional7}
\lim_{m \rightarrow \infty} \left| 
				{\langle \bar{\mathcal{R}}^m (u)( \bar{v}, \bar{w} ) \bar{v} \, , \, \bar{w}  \rangle}_{L^2 (\rho^m dx)} + 
				\int_{S_t} \rho_+ \nabla_{N_+} \Delta^{-1} \tr{(Dv_+)}^2 {|w_+^{\bot}|}^2 \, dS \, \right| 
				\leq  C {| w |}^2_{L^2 (\R^n \minus S_t) } \, .
\end{eqnarray}
This shows that the leading order term of the sectional curvature of $\G$ in the limit $\rho_- \rightarrow 0$
is given by the self--adjoint operator $\bar{\mathcal{R}}_0 (v)$ acting on $T_{u_+} \G$ 
represented in Lagrangian coordinates by
\begin{equation*}
\bar{\mathcal{R}}_0 (v_+) =   \left( - \, \rho_+ \nabla \H_+ \left( \nabla_{N_+} \Delta^{-1} \tr{(Dv_+)}^2    
						{(  \left. \cdot \,\, \right|_{\partial u_+ (\Omega_0) } )}^\bot	\right) \right) \circ u_+
\end{equation*}
and satisfying
\begin{equation}
\label{intRT2}
{\langle \bar{\mathcal{R}}_0 (v_+) \bar{w}_+ , \, \bar{w}_+   \rangle}_{L^2 (\rho_+ dy)}= - \int_{ \partial u_+ (\Omega_0) }
						 \nabla_{N_+} \Delta^{-1} \tr{(Dv_+)}^2 {|w_+^\bot|}^2 \rho_+ \, dS \, .
\end{equation}
From \eqref{pwaterwave} we see that $\Delta^{-1} \tr{(Dv_+)}^2$  is exactly $p^\star_{v_+,v_+}$ for the water wave problem 
so that \eqref{intRT2} is equivalent to the first integral in \eqref{intRT}.
Therefore we have shown that as $\rho_- \rightarrow 0$
the Kelvin--Helmotz instability for the vortex--sheet problem becomes the Raileigh--Taylor  instability, 
i.e. the leading order term of the sectional curvature of $\G$
is not definite in general and has a positive sign only provided that 
the normal gradient of the physical pressure (in absence of surface tension) is negative. 
We conclude with two observations:

\begin{itemize}

\item[(i)]
If we do not restrict our attention exclusively to the highest order term of the sectional curvature operator,
the above calculations show that
\begin{eqnarray*}
\lim_{m \rightarrow \infty} {
			\langle \bar{\mathcal{R}}^m (u)( \bar{v}, \bar{w} ) \bar{v} \, , \, \bar{w} \rangle}_{L^2 (\rho^m dy)} =
			-  \int_{S_t} \rho_+ \nabla_{N_+} \Delta^{-1} \tr{(Dv_+)}^2  \, {|w_+^\bot |}^2 \, dS 
\\
\\
			-  \int_{ \R^n \minus S_t } D^2 \Delta^{-1} \tr{(Dv_+)}^2 (w_+, w_+) \, \rho_+
			-  {|  \nabla \Delta^{-1} \tr( Dv_+ Dw_+ ) |}^2  \, \rho_+ \, dx \, .
\end{eqnarray*}
Since the second fundamental form of $\G^\star$ in the water wave problem is given by 
$\nabla p^\star_{v,w} = \nabla \Delta^{-1} \tr(Dv Dw)$ 
the above limit is exactly 
\begin{equation*}
\label{secww} 
\int_{ \R^n \minus S_t } \nabla p^\star_{v_+,v_+} \nabla p^\star_{w_+,w_+} \, \rho_+ \, dx  
					-   \int_{ \R^n \minus S_t } {|\nabla p^\star_{v_+, v_+} |}^2 \, \rho_+ \, dx  =
					{\langle \bar{\mathcal{R}}^\star (u_+)( \bar{v}_+, \bar{w}_+ ) \bar{v}_+  
					\, , \, \bar{w}_+ \rangle }_{L^2 (\rho_+ dy)}.
\end{equation*}
\item[(ii)]
From a standard argument we conclude that the full
curvature tensor of $\G$ converges to the curvature tensor of $\G^\star$
in the sense stated in \eqref{curvatureconv} 
and this completes the proof of proposition \ref{theocurvature}
$\, _\blacksquare$

\end{itemize}
\appendix

\section{Supporting material for proofs}

In this appendix we gather some technical results needed in the proofs presented and 
in the proof of theorem \ref{teoVS2} in appendix \ref{proofVS2}.

\begin{lem}
\label{circu}
Let $D_i$ (resp. $S_i$) be domains (resp. hypersurfaces) in $\R^n$ for $i=0,1$.
Let $\eta : D_0 \rightarrow D_1$ (resp. $\eta : S_0 \rightarrow S_1$) 
be an $H^l$--diffeomorphism for $l >\frac{n}{2} + 1$ (resp. $l >\frac{n+1}{2}$)
with ${| {(\det D \eta )}^{-1} |}_{L^\infty (D_1)} \leq a $ (resp. ${| {(\det D \eta)}^{-1} |}_{L^\infty (S_1)} \leq a $).  
Then the operator 
$T_\eta : f \rightarrow f \circ \eta$ is a bounded operator from $H^s (D_1)$ to $H^s (D_0)$ 
(resp. from $H^s (S_1)$ to $H^s (S_0)$) for any $s \in [0,l]$ and satisfies
\begin{equation}
\label{ineqcircu}
{ | f \circ \eta | }_{ H^s (D_0) }  \leq  C_0 { | f | }_{ H^s (D_1) } 
																								{ | \eta | }^s_{ H^l (D_0)}
\end{equation}
for some constant $C_0$ depending on $a$, $s$, $l$ and the domains $D_i$ (resp. the hypersurfaces $S_i$).
\end{lem}

\pr
The case $s=0$ follows immediately from the hypotheses.
Assume by induction that (\ref{ineqcircu}) holds for any integer $s$ such that $0 \leq s \leq k - 1 \leq l - 1$.
We prove the statement for $s = k$.
Write $D^k (f \circ \eta) = D^{k-1} (Df \circ \eta \, D\eta) = \sum_{j=0}^{k-1} D^j ( Df \circ \eta) D^{k-j} \eta$.
Let $r \geq 2$ be the integer such that $\frac{n}{2} - 1 \leq  l - r  < \frac{n}{2}$; 
observe that $D^i \eta \in L^\infty$ for $i \leq r - 1$ 
while it is not uniformly bounded in general for $i \geq r$ since $H^{l-i}$ does not embed in $L^{\infty}$.
According to this we split
\begin{equation*}
\sum_{j=0}^{k-1} D^j ( Df \circ \eta) D^{k-j} \eta =  \sum_{j=0}^{k-r} D^j ( Df \circ \eta) D^{k-j} \eta 
																											+ \sum_{j=k-r+1}^{k-1} D^j ( Df \circ \eta) D^{k-j} \eta
																											=: \Sigma_1 + \Sigma_2 \, .
\end{equation*}
In $\Sigma_2$ all derivatives on $\eta$ can be taken in $L^\infty$ and estimated 
through Sobolev's embedding:
\begin{eqnarray*}
{| \Sigma_2 |}_{L^2 (D_0)} & \leq & \sum_{j=k-r+1}^{k-1}  {| D^j ( Df \circ \eta) |}_{L^2 (D_0)} 
																											{|D^{k-j} \eta|}_{L^\infty (D_0)} \leq
																											C {| Df \circ \eta |}_{H^{k-1} (D_0)}
																											{|\eta|}_{H^l (D_0)}
\\
\\
& \leq &	C {| Df |}_{H^{k-1} (D_1)} {| \eta |}^{k-1}_{H^l (D_0)}	{|\eta|}_{H^l (D_0)}
					= C  {| f |}_{H^k (D_1)} {| \eta |}^k_{H^l (D_0)}   \, .
\end{eqnarray*}
The contribution of $\Sigma_1$ is estimated using H\"{o}lder's inequality and Sobolev's embeddings.
Since $l > \frac{n}{2} + 1$ and $k - j \geq r > 1$, 
we can choose $2 < p,q < \infty$ such that
\begin{equation*}
\frac{1}{p} + \frac{1}{q} = \frac{1}{2}  
											\hskip4pt  , \hskip10 pt 
											\frac{1}{q} > \frac{1}{2} - \frac{l - k + j}{n} 
											\hskip4pt  , \hskip10 pt 
											\frac{1}{p} > \frac{1}{2} - \frac{k - 1 - j}{n} \, .
\end{equation*}
Using H\"{o}lder's inequality and the embeddings
$H^{l-k+j} \subset H^{(\frac{1}{2} - \frac{1}{q})n } \subset L^q$,  
$H^{k-1-j} \subset H^{(\frac{1}{2} - \frac{1}{p})n } \subset L^p $ we get
\begin{eqnarray*}
{| \Sigma_1 |}_{L^2 (D_0)} \leq \sum_{j=0}^{k-r}  {| D^j ( Df \circ \eta) |}_{L^p (D_0)} 
																											{|D^{k-j} \eta|}_{L^q (D_0)} \leq
																											C {| D^j ( Df \circ \eta ) |}_{H^{k-1-j} (D_0)}
																											{|\eta|}_{H^l (D_0)}
													 \leq C  {| f |}_{H^k (D_1)} {| \eta |}^k_{H^l (D_0)}
\end{eqnarray*}
with $C$ depending on $a$, $k$ and the domains $D_0$, $D_1$. Therefore we proved 
\begin{equation*}
{ |D^k (f \circ \eta) | }_{L^2 (D_0)}  \leq  C  {| f |}_{H^k (D_1)} {| \eta |}^k_{H^l (D_0)}
\end{equation*}
From the inductive hypothesis the same inequality holds for lower order derivatives 
terms $D^i (f \circ \eta)$, for $0 \leq i \leq k-1$, replacing $k$ with $i$. 
Up to further increasing the value of $C$ depending on $a$ and the constants in Sobolev's embeddings,
we can sum this inequalities to obtain \eqref{circu} for $s=k$.
The case for non-integer $s$ follows by interpolation $_\Box$

\begin{lem}[{\bf About differential operators on $\Lambda_0$ \cite{shatah1}}]
\label{lemmaoperator}
Let $\Delta^{-1}$ and $\H$ denote respectively the inverse Laplacian with Dirichlet boundary condition
and the Harmonic extension\footnote{
$F = \Delta^{-1} f$ satisfies $\Delta F = f$ in $\Omega$ and $F = 0$ on $\partial \Omega$.
$G = \H g$ satisfies $\Delta G = 0$ in $\Omega$ and $G = g$ on $\partial \Omega$.
}
operators on a domain $\Omega$.
There exists a uniform constant $C > 0$ such that for every domain $\Omega$ 
with $\partial \Omega = S \in \Lambda_0$ (see definition \ref{defLambda}) the following is true:
\begin{eqnarray}
\label{estimaterestriction}
{| {\left. f  \right |}_{\partial \Omega}  |}_{H^{s} (\partial \Omega) } & \leq & 
															C {| f |}_{H^{s + \frac{1}{2}} (\Omega)}
																					\hskip8pt , \hskip8pt \for s > 0
\\
\label{estimateDH}
{| \Delta^{-1} |}_{L ( H^{s-1} (\Omega), H^{s+1} (\Omega ) )}  & + &
{| \H |}_{L ( H^{s + \frac{1}{2} } (\partial \Omega), H^{s+1} (\Omega ) )}  \leq C  \hskip8pt , \hskip8pt \for s \in (0,l-1] \, .
\end{eqnarray}
As a consequence the Dirichlet--to--Neumann operator satisfies\footnote{
In view of \eqref{estimateDH} $\mathcal{N}_0$ can be defined
for any $f \in H^s (\partial \Omega)$, $s \geq \frac{1}{2}$ in the weak form
$\langle \varphi , \mathcal{N}_0 ( f ) \rangle = \int_{\Omega} \nabla \H \varphi \nabla \H f $.
}
\begin{equation}
\label{estimateN0}
{| \mathcal{N}_0 |}_{L ( H^{s + \frac{1}{2} } (\partial \Omega), H^{s - \frac{1}{2} } (\partial \Omega ) )}  +
							{| \mathcal{N}_0^{-1} |}_{L ( \dot{H}^{s - \frac{1}{2} } (\partial \Omega), 
																						\dot{H}^{s + \frac{1}{2}} (\partial \Omega) )}
							\leq C  \,\,  ,
							\hskip6pt \for  s \in [0, l - 1]  \, ,
\end{equation}
where $\dot{H}^s$ denotes zero--mean $H^s$--functions.
In particular if $\mathcal{N}$ is the operator defined by
\eqref{N} 
then  $ {| \mathcal{N} |}_{L ( H^{s + \frac{1}{2}} (\partial \Omega), H^{s-\frac{1}{2}} (\partial \Omega ) )}  
\leq C (\rho_- + \rho_+)/ (\rho_- \rho_+)$ and 
\begin{equation}
\label{estimateN-1}
{ |  \mathcal{N}^{-1}  | }_{L ( \dot{H}^{s - \frac{1}{2} } 
															(\partial \Omega),  \dot{H}^{s + \frac{1}{2}} (\partial \Omega ) )} 
							\leq 2 C \rho_-	 \,\,  ,  \hskip6pt \for  s \in [0, l - 1]  \,\, 
							\mbox{and} \,\, \rho_- \leq \frac{\rho_+}{2 C^2} \, .
\end{equation}
\end{lem}

\pr
The proof of \eqref{estimaterestriction},  \eqref{estimateDH}, \eqref{estimateN0} 
and more detailed analysis of operators acting on $\partial \Omega$ 
(and in particular of the Dirichlet--to--Neumann operator) can be found in \cite[A.2]{shatah1}.
To obtain \eqref{estimateN-1} write $\mathcal{N}$ as
\begin{equation*}
\mathcal{N} =  
\left(   \frac{\mathcal{N}_+}{\rho_+}  \rho_- \mathcal{N}_-^{-1} +  I  \right) \frac{\mathcal{N}_-}{\rho_-}
=: (B + I) \frac{\mathcal{N}_-}{\rho_-}  \, .
\end{equation*}
Estimate \eqref{estimateN0} implies that for $\rho_- \leq \rho_+/(2 C^2)$ the linear operator $B$
maps $H^s (\partial \Omega)$ to itself with norm less or equal than $\frac{1}{2}$.
Thus $I + B$ is invertible and
$ \mathcal{N}^{-1} = \rho_- \mathcal{N}_-^{-1} \sum_{j=0}^\infty {(-1)}^j B^j $
so that
\begin{equation*}
{| \mathcal{N}^{-1} |}_{L (H^{s - \frac{1}{2}} (\partial \Omega), H^{s+\frac{1}{2}} (\partial \Omega))} 
										\leq  \rho_- C \sum_{j=0}^\infty 
										{| B |}^j_{L (H^{s - \frac{1}{2}} (\partial \Omega), H^{s+\frac{1}{2}} (\partial \Omega))} 
										\leq 2 C \rho_- \,\,\,  _\Box
\end{equation*}

\begin{lem}[{\bf Geometric Formulae \cite{shatah1}}]
Let $S$ be an hypersurface in $\R^n$ moved by the normal component of a vector field $v$.
Let $N, \k$ and $\Pi$ denote respectively its unit normal, mean curvature and second fundamental form. 
Then the following identities hold true:
\begin{eqnarray}
\label{D_tN}
\mathbf{D}_t N  & =	& - { \left[ (Dv)^\star \cdot N \right]}^\top
\\
\nonumber
\\
\label{D_tkappa}
\mathbf{D}_t \kappa  & = & - \Delta_{S_t} v^\bot - v^\bot {|\Pi|}^2 + \nabla_{v^\top} \kappa
\\
\nonumber
\\
\label{D_tPI}
\mathbf{D}_t^\top \Pi (\tau) & = &  - \D_\tau \left(  {( (Dv)^\star N)} ^\top \right) 
																	- \Pi \left( {( \nabla_\tau v )}^\top \right)
\\
\nonumber
\\
\label{formulaPIk}
- \Delta_S \Pi & = & - \D^2 \k  +  ( {|\Pi|}^2 I - \k \Pi ) \Pi
\end{eqnarray}
where $\mathcal{D}$ denotes the covariant derivative on $S$ and $\Delta_S := \tr \D^2$.
Furthermore there exists a uniform constant $C$ such that for any $S \in \Lambda_0$
\begin{equation}
\label{stimaNPI}
{|\Pi|}_{H^s (S)} + {|N|}_{H^{s+1} (S)} \leq C (1 + {|\kappa|}_{H^s(S)}) 
																			\hskip10pt \for l - \frac{5}{2} \leq s \leq l - 1 \, .
\end{equation}
\end{lem}
\begin{lem}[{\bf Commutator Estimates \cite{shatah1}}]
There exists a uniform constant $C$ such that for any $\partial \Omega = S \in \Lambda_0$ the following estimates hold:
\begin{eqnarray}
\label{D_tH}
{ \left|  [ \mathbf{D}_{t}, \H ]   \right|  }_{L (H^{s - \frac{1}{2}} (S) , H^{s} (S)) } 
																								& \leq & C {|v|}_{H^{l} (\Omega)}
																								\hskip10pt \for \frac{1}{2} < s \leq l
\\
\nonumber
\\
\label{D_tDelta-1}
{ \left| \left[ \mathbf{D}_{t}, \Delta^{-1} \right] \right| }_{ 
																								L(H^{s - 2} (\Omega), H^s (\Omega))}
																								& \leq & C {|v|}_{H^{l} (\Omega)}
												\hskip10pt \for 2 - l \leq s \leq l
\\
\nonumber
\\
\label{commest1}
{ \left| \left[ \mathbf{D}_{t}, \mathcal{N}_0 \right] \right| }_{
																								L(H^{s} (S), H^{s-1} (S))}
																								& \leq & C {|v|}_{H^{l} (\Omega)}
																								\hskip10pt \for 1 \leq s \leq l - \frac{1}{2}
\\
\nonumber
\\
\label{commest2}
{ \left| \left[ \mathbf{D}_t, \Delta_{S} \right] \right| }_{L (H^{s} (S), H^{s-2} (S) )}
																								& \leq & C {|v|}_{H^{l} (\Omega)}
																								\hskip10pt \for \frac{7}{2} - l < s \leq l - \frac{1}{2}  \, .
\end{eqnarray}
\end{lem}

\section{Proof of Theorem \ref{teoVS2}}
\label{proofVS2}
This section  is devoted to the proof of Theorem \ref{teoVS2} 
and consists essentially of material contained in \cite[sec. 4.3, 4.4]{shatah2}.
The only difference is that we claim and show independence of the energy estimates
on the densities of the two fluids.
Therefore, even though the proof is extremely similar to the one performed in \cite{shatah2}, 
we present it here for the reader's convenience.

\subsection{Estimates on the Lagrangian coordinate map}
We use the same notation in the original proof of theorem \ref{teoVS2}
letting $l := \frac{3}{2}k $. 
Working on the compact domain $\Omega_t^+$ and arguing as in the proof of proposition \ref{proVS1} 
(see section \ref{proofVS1})
we obtain the existence of a positive time $t_1$ and a constant $C_1$, 
only depending on $k,n$ and $\mu$ as in \eqref{t_0} such that
\begin{equation*}
\label{t_1}
{| u_+ (t, \cdot) - \id_{\Omega_0^+} | }_{H^{\frac{3}{2}k} (\Omega_0^+) } \leq C_1 t  
																			\hskip10pt \for t \in [0, \min\{t_0,t_1\}] \, .
\end{equation*}
This implies the estimate on the mean curvature\footnote{
This can be checked using the local coordinates constructed in \cite[appendix A]{shatah1}.
} 
\begin{equation}
\label{estklow}
{| \k_+ (t, \cdot) | }_{H^{\frac{3}{2} k - \frac{5}{2} } (S_t) } \leq C_2 t + 
																{| \k_+ (0, \cdot) | }_{H^{\frac{3}{2} k - \frac{5}{2} } (S_0) }	
 																														\hskip10pt \for t \in [0, \min\{t_0,t_1\}] \, .
\end{equation}
where the constant $C_2$ is only determined by $\mu$ and the set
$\Lambda_0$.
We conclude that there exists a time $t_2$ again determined only by $\mu$ and the set $\Lambda_0$ such that 
\begin{equation*}
S_t \in \Lambda_0 		\hskip8pt , \hskip8pt  \for t \in [ 0, \min\{t_0,t_2\}] \, .
\end{equation*}

\subsection{Evolution of the Energy}
The energy defined in \eqref{energyVS} is made of three terms. 
The first two involve the operator $\mathpzc{A}$ defined in \eqref{A} and are used to control the
irrotational part of the velocity and the mean curvature (hence the regularity of the evolving domain $S_t$);
the third part involves the vorticity $\o$ and is used to control the rotational part of $v$.
More explicitly $E = E_1 + E_2 + {|\o|}^2_{H^{\frac{3}{2}k - 1}}$ 
where, using \eqref{A}, the first two terms are given by
\begin{eqnarray}
\label{E_1}
E_1 & := & \frac{1}{2} \int_{\R^n \minus S_t} {| \mathpzc{A}^{ \frac{k}{2} } v |}^2 \rho \, dx 
=
\frac{1}{2} \int _{S_t}  v_+^\bot  {(- \Delta_{S_t} \bar{\mathcal{N}}) }^{k-1} (- \Delta_{S_t}) v_+^\bot \, dS
\\
\nonumber
\\
\label{E_2}
E_2 & := & \frac{1}{2} \int_{\R^n \minus S_t}  {| \mathpzc{A}^{ \frac{k}{2} - \frac{1}{2} } \nabla p_\k |}^2 \rho \, dx
= 
\frac{1}{2} \int_{S_t}   \k_+ \bar{\mathcal{N}} {(- \Delta_{S_t} \bar{\mathcal{N}}) }^{k-1} \k_+ \, dS
\end{eqnarray}
where
\begin{equation}
\label{barN}
\bar{\mathcal{N}} = 
\left( \frac{1}{\rho_+} \mathcal{N} \right) \mathcal{N}^{-1} \left( \frac{1}{\rho_-} \mathcal{N} \right) \, .
\end{equation}
It is clear from lemma \ref{lemmaoperator} that $\bar{\mathcal{N}}$ is a first--order self--adjoint operator
whose norm and inverse's norm do not depend on $\rho_-$.

\begin{pro}
There exists a polynomial 
$ 
Q (t) = Q \left( \, {|v (t, \cdot)|}_{ H^{\frac{3}{2}k } (\R^n \minus S_t) } , \, 
{|\kappa (t, \cdot)|}_{ H^{\frac{3}{2}k - 1} (S_t) }  \right) 
$
with positive coefficients depending on the set $\Lambda_0$ and independent of the density $\rho_-$ 
such that
\begin{equation}
\label{evolE}
\left|  \frac{d}{dt} (E - E_{\ex}) \right| \leq Q
\end{equation}
where the extra energy term $E_{\ex}$, due to the Kelvin--Helmotz instability, is given by
\begin{equation}
\label{E_ex}
E_{\ex} = \frac{\rho_+}{ 2 (\rho_+ + \rho_-) } 
					\int_{S_t} \nabla_{v_+^\top} \kappa_+ \cdot \bar{\mathcal{N}} 
					{(- \Delta_{S_t} \bar{\mathcal{N}}) }^{k-2} \nabla_{v_+^\top} \kappa_+  \, dS
					-  \frac{\rho_-}{ 2 (\rho_+ + \rho_-) } 
					\int_{S_t} \nabla_{v_-^\top} \kappa_+ \cdot \bar{\mathcal{N}} 
					{(- \Delta_{S_t} \bar{\mathcal{N}}) }^{k-2} \nabla_{v_-^\top} \kappa_+
					\, dS \, .				
\end{equation}
\end{pro}

\pr
Throughout the proof we denote by $Q$ any generic polynomial satisfying the properties in the statement.

\nl
{\it Evolution of $E_1$}:
This is the hardest term to deal with and is the one where the extra energy term $E_{\ex}$ appears.
We are going to show 
\begin{equation}
\label{evolE_1}
\left| \frac{d}{dt} \left( E_1 - E_{\ex} \right)
					+
					\int_{S_t}   
					v_+^\bot  {(- \Delta_{S_t} \bar{\mathcal{N}}) }^k  \k_+  \, dS
					\right| \leq Q	 \, .
\end{equation}
From definition \eqref{barN} and \eqref{commest1} we have 
\begin{equation*}
{ \left| \left[ \mathbf{D}_t, \bar{\mathcal{N}} \right] \right| }_{L (H^s (S_t), H^{s-1} (S_t) )}
																								\leq C {|v|}_{H^{\frac{3}{2} k} (\Omega_t)}
																								\hskip10pt \for  \frac{1}{2} \leq s \leq \frac{3}{2} k - \frac{1}{2}  \, .
\end{equation*}
Therefore, using also \eqref{commest2}, 
we can commute $\mathbf{D}_{t_+}$ with the operators appearing in \eqref{E_1} to get
\begin{equation*}
\left| \frac{d}{dt} \frac{1}{2} \int_{S_t}   
					v_+^\bot {(- \Delta_{S_t} \bar{\mathcal{N}}) }^{k-1} (- \Delta_{S_t} )  v_+^\bot \, dS  
					-
					\int_{S_t}   
					v_+^\bot {(- \Delta_{S_t} \bar{\mathcal{N}}) }^{k-1} (- \Delta_{S_t} ) \mathbf{D}_{t_+} v_+^\bot \, dS
					\right| \leq Q	\, .				
\end{equation*}
Using \eqref{D_tN} to express $\mathbf{D}_{t_+} N_+$, \eqref{Ematerial}, \eqref{p_vw} and \eqref{Sp}
together with \eqref{barN} we have
\begin{eqnarray*}
\mathbf{D}_{t_+} v_+^\bot & = &  \mathbf{D}_{t_+} v_+ \cdot N_+  +  v_+  \mathbf{D}_{t_+} N_+
													=  - \nabla_{N_+} p_{v,v}^+  -  \nabla_{N_+} p_\k^+  -  \nabla_{v_+^\top} v_+ \cdot N_+
													\\
													\\
													& = &  - \frac{1}{\rho_+} \mathcal{N}_+ p^S_{v,v}  +  \nabla_{N_+} \Delta_+^{-1} \tr {(Dv)}^2
													- \bar{\mathcal{N}} \k_+  -  \nabla_{v_+^\top} v_+^\bot  +  \Pi (v_+^\top, v_+^\top)
\end{eqnarray*}
From lemma \ref{lemmaoperator} and trace--estimates
$
{ | \nabla_{N_+} \Delta_+^{-1} \tr {(Dv)}^2  | }_{H^{ \frac{3}{2}k - \frac{1}{2} } (S_t) }  \leq Q
$ 
so that this term is lower order and
\begin{equation}
\left| \frac{d}{dt} E_1
					-
					\int_{S_t}   
					v_+^\bot {(- \Delta_{S_t} \bar{\mathcal{N}}) }^{k-1} (- \Delta_{S_t} )
					\left( - \frac{1}{\rho_+} \mathcal{N}_+ p^S_{v,v}   - \bar{\mathcal{N}} \k_+  
								 -  \nabla_{v_+^\top} v_+^\bot  +  \Pi (v_+^\top, v_+^\top) \right)
					\, dS  
					\right| \leq Q	\, .
\label{estE_11}					
\end{equation}
Equation \eqref{p_vw} for $p^S_{v,v}$ gives
\begin{eqnarray*}
- \frac{1}{\rho_\pm} \mathcal{N}_+ p^S_{v,w} & = & 
						- \frac{1}{\rho_\pm} \mathcal{N}_+ \mathcal{N}^{-1} \left\{ 2 \nabla_{v_+^\top - v_-^\top} v_+^\bot
				    - \Pi_+ ( v_+^\top, v_+^\top )
				    - \Pi_- ( v_-^\top, v_-^\top ) 
				    - \nabla_{N_+} \Delta_+^{-1}  \tr {(Dv)}^2  - \nabla_{N_-} \Delta_-^{-1}  \tr {(Dv)}^2
				    \right\} \, .
\end{eqnarray*}
Since $\mathcal{N}_+ \mathcal{N}^{-1}$ is an operator of order zero the terms 
$\nabla_{N_\pm} \Delta_{\pm}^{-1}  \tr {(Dv)}^2$ can be treated as before. 
From Lemma 4.6 in \cite{shatah2} (to which we refer for the proof)
\begin{equation*}
{ \left|  {( - \Delta_{S_t} )}^\frac{1}{2} - \mathcal{N}_\pm \right|  }_{L (H^s (S_t)) }  \leq
	 						    C \left( 1 +  {| \k (t, \cdot) |}_{H^{\frac{3}{2} k - \frac{3}{2} (S_t)} } \right)
									\hskip10pt \for  \frac{1}{2} - \frac{3}{2}k \leq s \leq \frac{3}{2}k - \frac{1}{2} \, ;
\end{equation*}
this and the definition \eqref{N} of $\mathcal{N}$ yield
\begin{equation*}
{ \left|  \mathcal{N}_+ \mathcal{N}^{-1}  -  \frac{ \rho_+ \rho_-}{\rho_+  +  \rho_-} \right|  
								}_{L \left( H^{\frac{3}{2} k - \frac{3}{2} } (S_t) , H^{\frac{3}{2} k - \frac{1}{2} } (S_t) \right) }  
								\leq  Q   \, .
\end{equation*}
Together with \eqref{stimaNPI} this gives
\begin{equation*}
{ \left|  - \frac{1}{\rho_+} \mathcal{N}_+ p^S_{v,v} -  
						\frac{ \rho_- }{\rho_+  +  \rho_-} \left(  2 \nabla_{v_+^\top - v_-^\top} v_+^\bot
				    - \Pi_+ ( v_+^\top, v_+^\top )  - \Pi_- ( v_-^\top, v_-^\top ) 
				    \right)
								\right|  }_{H^{\frac{3}{2} k - \frac{1}{2} } (S_t) }  
								\leq  Q
\end{equation*}
so that \eqref{estE_11} becomes
\begin{eqnarray*}
\left| \frac{d}{dt} E_1  \right.
					& - &
					\int_{S_t}   
					v_+^\bot { (- \Delta_{S_t} \bar{\mathcal{N}}) }^{k-1} (- \Delta_{S_t} )
					\left[   - \frac{ \rho_- }{\rho_+  +  \rho_-} \Pi_- ( v_-^\top, v_-^\top ) 
					+ \frac{ \rho_+ }{\rho_+  +  \rho_-} \Pi_+ ( v_+^\top, v_+^\top )
					\right.
					\\
					\\
					& - & 
					\left.
					\left.
					\bar{\mathcal{N}} \k_+  
				  + \nabla v_+^\bot  \left(  \frac{\rho_- - \rho_+ }{\rho_+ + \rho_-} v_+^\top
				  - \frac{2 \rho_-}{\rho_+ + \rho_-} v_-^\top  \right)  
				  \right]
					\, dS
					\right| \leq Q	 \, .
\end{eqnarray*}
We now claim that the last two terms in the above integral are lower order.
To see this, consider flows $\Phi_\pm (\tau, \cdot)$ on $\Omega_t^+$ 
generated by $\H_+ v_\pm^\top$ 
and apply \eqref{commest1} and \eqref{commest2} to $\mathbf{D}_\tau$
to move outside the tangential derivatives $\nabla_{v_\pm^\top}$:
\begin{eqnarray*}
& & \left| \int_{S_t} v_+^\bot { (- \Delta_{S_t} \bar{\mathcal{N}}) }^{k-1} (- \Delta_{S_t} )
					\nabla v_+^\bot  \left(  \frac{\rho_- - \rho_+ }{\rho_+ + \rho_-} v_+^\top
				  - \frac{2 \rho_-}{\rho_+ + \rho_-} v_-^\top  \right)  
					\, dS
					\right.
					\\
					\\
					& - &
					\frac{\rho_- - \rho_+ }{2 (\rho_+ + \rho_-)} \int_{S_t} \nabla_{v_+^\top} 
					{\left| { (- \Delta_{S_t} \bar{\mathcal{N}}) }^{\frac{k-1}{2}} {(- \Delta_{S_t} )}^{\frac{1}{2} } 
					v_+^\bot  \right|}^2 
					\, dS
				  \\
				  \\
				  & + &
				  \left. 
				  \frac{\rho_-}{\rho_+ + \rho_-} \int_{S_t} \nabla_{v_-^\top} 
				  {\left| { (- \Delta_{S_t} \bar{\mathcal{N}}) }^{\frac{k-1}{2}} {(- \Delta_{S_t} )}^{\frac{1}{2} } 
					v_+^\bot  \right|}^2
					\, dS
					\right| \leq Q	 \, ;
\end{eqnarray*}
then integrate by parts in these last two integrals
estimating 
$D v_\pm^\top$ in $L^\infty (S_t)$ 
and the remaining $\frac{3}{2} k - \frac{1}{2}$ derivatives on $v_+^\bot$ in $L^2 (S_t)$.

\nl
For the terms involving the second fundamental form, \eqref{formulaPIk} gives
\begin{equation*}
\left| \Delta_{S_t} ( \Pi_\pm (v_\pm^\top, v_\pm^\top) - \D^2 k_\pm (v_\pm^\top, v_\pm^\top)   \right|_{
								H^{\frac{3}{2}k - \frac{5}{2}} (S_t) }  \leq  Q   \,  .
\end{equation*}
Since 
$ \D^2 \k_\pm = \nabla_{v_\pm^\top} \nabla_{v_\pm^\top} \k_\pm - \D_{v_\pm^\top} v_\pm^\top \cdot \nabla \k_\pm$
and the last term in this sum is lower order, we get
\begin{eqnarray*}
\left|  \frac{d}{dt} E_1 - 
					\int_{S_t}   
					v_+^\bot { (- \Delta_{S_t} \bar{\mathcal{N}}) }^{k-1}  
				  \left[ \frac{ \rho_- }{\rho_+  +  \rho_-}  \nabla_{v_-^\top} \nabla_{v_-^\top} \k_-
					- 
					\frac{ \rho_+ }{\rho_+  +  \rho_-} \nabla_{v_+^\top} \nabla_{v_+^\top} \k_+
					+  \Delta_{S_t}  \bar{\mathcal{N}} \k_+  
				  \right]
					\, dS  
					\right| \leq  Q  \, .
\end{eqnarray*}
Using the same previous argument we can commute one of the factors $\nabla_{v_\pm^\top}$
and move it outside to obtain
\begin{eqnarray}
\left|  \frac{d}{dt} E_1 \right.
					& + & 
					\frac{ \rho_- }{\rho_+  +  \rho_-} 
					\int_{S_t}   
					\nabla_{v_-^\top}  ( - \Delta_{S_t} v_+^\bot ) \bar{\mathcal{N}} 
					{ (- \Delta_{S_t} \bar{\mathcal{N}}) }^{k-2}  \nabla_{v_-^\top} \k_- \, dS
					\label{estE_12}	
					\\
					\nonumber
					\\
					& - &
					\frac{ \rho_- }{\rho_+  +  \rho_-} 
					\int_{S_t}   
					\nabla_{v_+^\top}  ( - \Delta_{S_t} v_+^\bot ) \bar{\mathcal{N}} 
					{ (- \Delta_{S_t} \bar{\mathcal{N}}) }^{k-2}  \nabla_{v_+^\top} \k_+ \, dS
					\label{estE_13}					
					\\
					\nonumber
					\\
					\nonumber			
					& + & \left. \int_{S_t}  v_+^\bot { (- \Delta_{S_t} \bar{\mathcal{N}}) }^k \k_+  
					\, dS  
					\right| \leq  Q \, .
\end{eqnarray}
Now, thanks to identity \eqref{D_tkappa}
\begin{equation*}
{ \left| - \Delta_{S_t} v_+^\bot  - \mathbf{D}_{t_+} \k_+  \right| }_{H^{\frac{3}{2}k - 2} (S_t) }  \leq  Q
\end{equation*}
so that we can substitute $\mathbf{D}_{t_+} \k_+$  to  $- \Delta_{S_t} v_+^\bot$ in \eqref{estE_12} and \eqref{estE_13}.
The usual commutator estimates imply
\begin{eqnarray*}
\left|  \frac{d}{dt} E_{\ex}  \right. 
					& + & \frac{ \rho_- }{\rho_+  +  \rho_-}
					\int_{S_t}  
					\nabla_{v_-^\top}  \mathbf{D}_{t_+} \k_+  \bar{\mathcal{N}} 
					{ (- \Delta_{S_t} \bar{\mathcal{N}}) }^{k-2}  \nabla_{v_-^\top} \k_- \, dS
					\\
					\nonumber
					\\
					& - &
					\left.
					\frac{ \rho_+ }{\rho_+  +  \rho_-} 
					\int_{S_t}   
					\nabla_{v_+^\top}  \mathbf{D}_{t_+} \k_+ \bar{\mathcal{N}} 
					{ (- \Delta_{S_t} \bar{\mathcal{N}}) }^{k-2}  \nabla_{v_+^\top} \k_+ \, dS
					\right| 
					\leq Q
\end{eqnarray*}
and \eqref{evolE_1} follows.

\nl
{\it Evolution of $E_2$}:  
As before commutator estimates \eqref{commest1} and \eqref{commest2} give
\begin{equation*}
\left| \frac{d}{dt} E_2
					-
					\int_{S_t}   
					\k_+ \bar{\mathcal{N}} {(- \Delta_{S_t} \bar{\mathcal{N}}) }^{k-1} \mathbf{D}_{t_+} \k_+ \, dS
					\right| \leq Q					
\end{equation*}
and in view of \eqref{D_tkappa} and \eqref{stimaNPI} we obtain
\begin{eqnarray*}
& & \left| \frac{d}{dt} E_2
					-
					\int_{S_t}   
					\k_+ \bar{\mathcal{N}} {(- \Delta_{S_t} \bar{\mathcal{N}}) }^{k-1} (- \Delta_{S_t} ) v_+^\bot \, dS  
					\right| 
					\leq Q 
					+  
					\left|  \int_{S_t}   
					\k_+ \bar{\mathcal{N}} {(- \Delta_{S_t} \bar{\mathcal{N}}) }^{k-1} \nabla_{v_+^\top} \k_+  \, dS \right| 					
\end{eqnarray*}
The same commutation argument previously adopted shows that
\begin{eqnarray*}
\left|  \int_{S_t}   
					\k_+ \bar{\mathcal{N}} {(- \Delta_{S_t} \bar{\mathcal{N}}) }^{k-1} \nabla_{v_+^\top} \k_+  \, dS
					- \frac{1}{2}
					\int_{S_t}   
					\nabla_{v_+^\top} { \left| \bar{\mathcal{N}}^\frac{1}{2} 
														{ (- \Delta_{S_t} \bar{\mathcal{N}}) }^{ \frac{k-1}{2} } \k_+   \right|}^2 \, dS
					\right| 	\leq Q  \, .
\end{eqnarray*}
Integrating by parts and estimating $ D v_+^\top $ in $L^\infty (S_t)$ 
and the remaining $\frac{3}{2} k - 1$ derivatives on $\k_+$ in $L^2$ shows that this last integral is bounded by $Q$.
Finally use the self--adjointness of $\bar{\mathcal{N}}$ and $\Delta_{S_t}$ to obtain
\begin{equation}
\label{evolE_2}
\left| \frac{d}{dt}  E_2
					-
					\int_{S_t}   
					v_+^\bot  {(- \Delta_{S_t} \bar{\mathcal{N}}) }^k  \k_+  \, dS  
					\right| \leq Q	 \, .
\end{equation}

\nl
{\it Evolution of the vorticity $\o = Dv - {(Dv)}^\star$}: 
Commuting repeatedly $\mathbf{D}_t$ with $D$ and using the identity
\begin{equation*}
\mathbf{D}_t \o = D \mathbf{D}_t v - {(Dv)}^2 - {(D \mathbf{D}_t v)}^\star  +  {({(Dv)}^\star)}^2
								= {({(Dv)}^\star)}^2 - {(Dv)}^2 = - \o Dv  - {(Dv)}^\star \o 
\end{equation*}
we have
\begin{eqnarray}
\nonumber
& & \frac{d}{dt}  \int_{\R^n \minus S_t} {|D^{\frac{3}{2} k - 1} \o |}^2 \, dx
						= \int_{\R^n \minus S_t} \mathbf{D}_t {|D^{\frac{3}{2} k - 1} \o |}^2 \, dx
						\leq  
						C  {|v (t, \cdot)|}_{ H^{\frac{3}{2}k} (\R^n \minus S_t)  } 
										  {|\o (t, \cdot)|}_{ H^{\frac{3}{2}k - 1} (\R^n \minus S_t) }  
										  \leq Q  \, .
\label{evolvorticity}
\end{eqnarray}
Summing up \eqref{evolE_1}, \eqref{evolE_2} and \eqref{evolvorticity} we get the desired cancellations
giving \eqref{evolE} 
$_\Box$

\subsection{The Energy Inequality}
Integrating in time \eqref{evolE} gives
\begin{equation}
E (t) - E(0) - E_{\ex} (t) + E_{\ex} (0) \leq  
		\int_0^t Q \left( {|v (s, \cdot)|}_{ H^{\frac{3}{2}k (\R^n \minus S_s)} }  , \, 
											{|\kappa (s, \cdot)|}_{ H^{\frac{3}{2}k - 1} (S_s) }   \right) 
		\, ds
\label{Einequality}
\end{equation}
for any $0 \leq t \leq \min \{ t_0, t_2 \}$. 
We can estimate the extra energy term \eqref{E_ex} by
\begin{eqnarray*}
| E_{\ex} (t) | & \leq &
		\frac{1}{2} \int_{S_t} {\left| \bar{\mathcal{N}}^\frac{1}{2} {(\Delta_{S_t} \bar{\mathcal{N}} )}^{\frac{k}{2} - 1} 
																															\nabla_{v_+^\top} \k_+  \right| }^2 \, dS
	+ \frac{1}{2} \int_{S_t} {\left| \bar{\mathcal{N}}^\frac{1}{2} {(\Delta_{S_t} \bar{\mathcal{N}} )}^{\frac{k}{2} - 1} 
																															\nabla_{v_-^\top} \k_+  \right| }^2 \, dS
								\\
								\\																															
								& \leq & 
								C {|v_\pm^\top \cdot \nabla \k_+|}_{ H^{\frac{3}{2}k - \frac{5}{2}} (S_t) }
								\leq
								C {|v (t, \cdot)|}_{ H^{\frac{3}{2}k - \frac{5}{8} 
								} (\R^n \minus S_t) }  \, 
											{|\kappa (t, \cdot)|}_{ H^{\frac{3}{2}k - \frac{3}{2}} (S_t) }   
\end{eqnarray*}
where the positive constant $C$ depends only on the set $\Lambda_0$.
Interpolating $v$ between $H^{ \frac{3}{2} k - \frac{3}{2} }$ and $H^{\frac{3}{2} k}$
and $\k$ between $H^{ \frac{3}{2} k - \frac{5}{2} }$ and $H^{ \frac{3}{2} k - 1}$
yields
\begin{equation*}
| E_{\ex} | \leq \frac{1}{2} E  +  C_1 \left( 1 + {|v|}^m_{ H^{\frac{3}{2}k - \frac{3}{2} } (\R^n \minus S_t)} 
							\right)
\end{equation*}
for some integer $m$ where the constant $C_1$, 
which includes ${|\k|}_{H^{\frac{3}{2}k - \frac{5}{2}} }$, 
depends ultimately only on $E_0$ and $\Lambda_0$ in view of \eqref{estklow}.
Using Euler equations \eqref{Ematerial} and lemma \ref{estimatep} to estimate the pressure, we have
\begin{equation*}
{ |\mathbf{D}_t v | }_{H^{\frac{3}{2}k - \frac{3}{2}} (\R^n \minus S_t) } =  
										\frac{1}{\rho_+} { |\nabla p_+ | }_{H^{\frac{3}{2}k - \frac{3}{2}} (\Omega_t^+) }
										+
										\frac{1}{\rho_-} { |\nabla p_- | }_{H^{\frac{3}{2}k - \frac{3}{2}} (\Omega_t^-) }
										\leq Q \,  .
\end{equation*}
We can then use the Lagrangian coordinate map to estimate
\begin{equation*}
\left| {| v (t, \cdot) |}_{H^{\frac{3}{2}k - \frac{3}{2}} (\R^n \minus S_t) }^m -
				{| v (0, \cdot) |}_{H^{\frac{3}{2}k - \frac{3}{2}} (\R^n \minus S_0) }^m \right| \leq \int_0^t Q(s) \, ds
\end{equation*}
and obtain
\begin{equation*}
| E_{\ex} | \leq \frac{1}{2} E  +  C_1 \left( 1 + {|v (0,\cdot)|}^m_{ 
							H^{\frac{3}{2}k - \frac{3}{2} } (\R^n \minus S_0)} \right)
							+ \int_0^t Q(s) \, ds
							\leq
							\frac{1}{2} E  +  C_2 + \int_0^t Q(s) \, ds
\end{equation*}
where $C_2$ is determined by $E_0$, the set $\Lambda_0$ and 
${|v (0,\cdot)|}_{H^{\frac{3}{2}k - \frac{3}{2} } (\R^n \minus S_0)}$.
Inserting this last inequality in \eqref{Einequality} we finally obtain
\begin{equation*}
E (S_t, v(t,\cdot)) \leq 3 E (S_0, v(0,\cdot)) + C_2 + \int_0^t Q(s) \, ds
\end{equation*}
for some $C_2$ as above.
Taking $\mu$ in \eqref{t_0} large enough compared to the initial data concludes the proof of theorem \ref{teoVS2}
$_\blacksquare$

\addcontentsline{toc}{section}{Bibliography}
\bibliographystyle{plain}
\bibliography{vortex_sheet_1}

\begin{thebibliography}{10}

\bibitem{Ambrose}
D.M. Ambrose.
\newblock Well-posedness of vortex sheets with surface tension.
\newblock {\em SIAMJ. Math. Anal.}, 35(1):211--244, 2003.

\bibitem{AmbroseMasmoudi2}
D.M. Ambrose and N.~Masmoudi.
\newblock Well-posedness of 3-d vortex sheets with surface tension.
\newblock {\em {C}ommun. Math. Sci.}, 5(2):391--430, 2007.

\bibitem{Arnold66}
V.I. Arnold.
\newblock {S}ur la géométrie differentielle des groups de {L}ie de dimension
  infinie et ses application à l'hydrodynamique des fluids parfait.
\newblock {\em {A}nn. Inst. Fourier (Grenoble)}, 16(1):319--361, {J}uly 1966.

\bibitem{Beale}
J.T. {B}eale{,} T.Y.~Hou and J.S. Lowengrub.
\newblock Growth rates for the linearized motion of fluid interfaces away from
  equilibrium.
\newblock {\em Comm. on Pure Appl. Math.}, 46(9):1269--1301, 1993.

\bibitem{Brenier}
Y.~Brenier.
\newblock Minimal geodesics on groups of volume-preserving maps and generalized
  solutions of the {E}uler equations.
\newblock {\em Comm. Pure Appl. Math.}, 52(4):411--452, {J}uly 1999.

\bibitem{CoutShko1}
A.~{C}heng{,} D.~Coutand and S.~Shkoller.
\newblock On the motion of {V}ortex {S}heets with surface tension.
\newblock {\em Comm. on Pure and Appl. Math.}, 61(12):1715--1752, December
  2008.

\bibitem{coutshko3}
A.~{C}heng{,} D.~Coutand and S.~Shkoller.
\newblock On the limit as the density ratio tends to zero for two perfect
  incompressible 3-{D} fluids separated by a surface of discontinuity.
\newblock {\em Preprint}, August 2009.
\newblock Available on the web at
  http://www.math.ucdavis.edu/~shkoller/pub/pub.html.

\bibitem{ChriLind}
D.~Christodoulou and H.~Lindblad.
\newblock On the motion of the free surface of a liquid.
\newblock {\em Comm. on Pure Appl. Math.}, 53(12):1536--1602, 2000.

\bibitem{CoutShko2}
D.~Coutand and S.~Shkoller.
\newblock Well-posedness of the free-surface incompressible {E}uler equations
  with or without surface tension.
\newblock {\em Journal of. Amer. Math. Soc.}, 20(3):829--930, July 2007.

\bibitem{Ebin1}
G.~Ebin.
\newblock {T}he equations of motion of a perfect fluid with free boundary are
  not well posed.
\newblock {\em {C}omm. Partial Differential Equations}, 12(10):1175--1201,
  1987.

\bibitem{Ebin2}
G.~Ebin.
\newblock {I}ll{-}posedness of the {R}aileigh{-T}aylor and {K}elvin{-H}elmotz
  problems for incompressible fluids.
\newblock {\em {C}omm. Partial Differential Equations}, 13(10):1265--1295,
  1988.

\bibitem{EbinMarsden}
G.~Ebin and J.~Marsden.
\newblock {G}roups of {D}iffeomorphisms and the {M}otion of an {I}ncompressible
  {F}luid.
\newblock {\em {A}nnals of math}, 92(1):102--163, {J}uly 1970.

\bibitem{Lindblad}
H.~Lindblad.
\newblock Well-posedness for the motion of an incompressible liquid with free
  surface boundary.
\newblock {\em {A}nn. of Math.}, 162(1):109--194, 2005.

\bibitem{shatah1}
J.~Shatah and C.~Zeng.
\newblock {G}eometry and a priori estimates for free boundary problems of the
  {E}uler equation.
\newblock {\em {C}omm. on pure and appl. math}, 61(5):698--744, {M}ay 2008.

\bibitem{shatah2}
J.~Shatah and C.~Zeng.
\newblock A priori estimates for fluid interface problems.
\newblock {\em {C}omm. on pure and appl. math}, 61(6):848--876, {J}une 2008.

\bibitem{shatah3}
J.~Shatah and C.~Zeng.
\newblock Local well-posedness for the fluid interface problem.
\newblock {\em Preprint}, 2009.

\bibitem{Shnirelman}
A.~Shnirelman.
\newblock The geometry of the group of diffeomorphisms and the dynamics of an
  ideal incompressible fluid.
\newblock {\em {M}at. {S}b. ({N}.{S}.)}, 128(1):82--109, 144, {J}uly 1985.

\bibitem{Wu1}
S.~Wu.
\newblock {W}ell-posedness in {S}obolev spaces of the full water wave problem
  in 2-d.
\newblock {\em {I}nvent. math}, 130(1):39--72, 1997.

\bibitem{Wu2}
S.~Wu.
\newblock {W}ell-posedness in {S}obolev spaces of the full water wave problem
  in 3-d.
\newblock {\em {J}. {A}mer. {M}ath {S}ociety math}, 12(2):445--495, 1999.

\end{thebibliography}
\end{document}